\documentclass{amsproc}

\date{ } \vskip 0.5cm \baselineskip 18pt
\newtheorem{define}{{\bf Definition}}

 \newcommand{\F}[2]{\frac{#1}{#2}}
 
\newcommand{\E}[0]{\varepsilon}
\newcommand{\lam}[0]{\lambda}
 \newcommand{\BE}[1]{\begin{array}{#1}}
\newcommand{\EE}[0]{\end{array}}
\makeatletter \@addtoreset{equation}{section}
\makeatother

\newcommand{\usection}[1]{section}
\newlength{\picwidth} \newlength{\picwidthadj}
\newcommand{\unipic}[2]{\vspace{#2}
\newcommand \epf{\rule{0.2cm}{0.2cm}\hfill \bigskip \par}
\setlength{\picwidth}{\the\textwidth} \setlength{\picwidthadj}{0pt}
\addtolength{\picwidth}{\picwidthadj} \special{dvitops:import #1
\the\picwidth #2}
\bigskip }

\newtheorem{theorem}{Theorem}[section]
\newtheorem{lemma}[theorem]{Lemma}

\theoremstyle{definition}

\theoremstyle{remark}

\numberwithin{equation}{section}

\begin{document}

\title{Compensated Compactness Method on Non-isentropic
Polytropic Gas Flow}

\author{Yun-guang Lu}
\address{K.K.Chen Institute for Advanced Studies, Hangzhou Normal University, P.R. China, email: ylu2005@ustc.edu.cn}



\begin{abstract}
In this paper, we are concerned with a model of polytropic gas flow, which consists the mass equation, the momentum equation and a varying entropy equation. First, a new technique, to set up a relation between the Riemann invariants of the isentropic system and the entropy variable $s$, coupled with the maximum principle, is introduced to obtain the a-priori $L^{\infty}$ estimates for the viscosity-flux approximation solutions. Second, the convergence framework from the compensated compactness theory on the system of isentropic  gas dynamics is applied to prove the pointwise convergence of the approximation solutions and the global existence of bounded entropy solutions for the Cauchy problem of the system with bounded initial data.  Finally, as a by-product, we obtain a non-classical bounded generalized solution $(\rho,u,s),$  of the original non-isentropic polytropic gas flow, which satisfies the mass equation and the momentum equation, the entropy equation with an extra nonnegative measure in the sense of distributions.
\end{abstract}
 \maketitle

 \section{ \bf Introduction}
The Euler system for compressible polytropic gas flow in one-space dimension
is the following system of three conservation laws
\begin{equation}\left\{
\label{1.1}
\begin{array}{l}
\rho_{t}+( \rho u)_{x}=0,         \\\\
       ( \rho u)_t+( \rho u^2+ k \rho^{\gamma} e^{\F{s}{c}})_x=0, \\\\
       (\F{1}{2} \rho u^2+ \F{k}{\gamma-1} \rho^{\gamma} e^{\F{s}{c}} )_t+
       (u( \F{1}{2} \rho u^2+ \F{k \gamma}{\gamma-1} \rho^{\gamma} e^{\F{s}{c}}))_x=0,
\end{array}\right.
\end{equation}
where $\rho$ is the density of gas, $u$ the velocity, $s$ the entropy and
$ \gamma > 1$ corresponds to the adiabatic exponent, $c> 0$ is the specific heat at constant volume, $ k>0 $ can be any constant under scaling. Without any loss
of generality, we may choose $ k= \F{ \theta^{2}}{\gamma}, \theta=\F{\gamma-1}{2} $ and $ c= \F{1}{2}$ for the simplicity.

For a smooth solution, the third equation in (\ref{1.1}) is equivalent to
\begin{equation}
\label{1.2}
      (\rho s)_{t}+( \rho u s)_{x}=0,
\end{equation}
which, united with the first and the second equations in (\ref{1.1}), forms
the following system of three conservation laws
\begin{equation}\left\{
\label{1.3}
\begin{array}{l}
\rho_{t}+( \rho u)_{x}=0,         \\\\
       ( \rho u)_t+( \rho u^2+ \F{ \theta^{2}}{\gamma} \rho^{\gamma} e^{2s})_x=0, \\\\
       (\rho s)_{t}+( \rho u s)_{x}=0.
\end{array}\right.
\end{equation}
System (\ref{1.3}) was derived by relaxation from an isentropic two-phase mixture, and
its global weak solutions were studied in \cite{1}, under the assumption of the uniform
boundedness of the viscosity solutions. In this paper, we shall give a rigorous proof of the uniform bound of the viscosity solutions and obtain  the global bounded entropy solutions of the Cauchy problem (\ref{1.3})  with the following bounded initial data
\begin{equation}
\label{1.4} (\rho,u,s)|_{t=0}=(\rho_{0}(x),u_{0}(x),s_{0}(x)), \hspace{0.3cm} \rho_{0}(x) \ge 0.
\end{equation}
When $ s$ is a constant,  both (\ref{1.1}) and
(\ref{1.3}) are reduced to the following isentropic gas dynamics system of two conservation laws
\begin{equation}\left\{
\label{1.5}
\begin{array}{l}
      \rho_{t}+(\rho u)_{x}=0,        \\\\
       ( \rho u)_t+( \rho u^2+ \F{ \theta^{2}}{\gamma} \rho^{\gamma} e^{2s} )_x =0.
\end{array}\right.
\end{equation}
Numerous papers deal with the analysis of weak solutions of the Cauchy problem (\ref{1.5}). The first existence theorem for large initial data of locally
finite total variation was proved in \cite{2} for $\gamma=1$ and in
\cite{3} for $ \gamma \in (1, 1+ \delta)$ in Lagrangian coordinates,
where $ \delta$ is small. The Glimm scheme \cite{4} was used  in these papers.

The ideas of compensated compactness developed in \cite{5,6}
were used  in \cite{7} to
 established a global existence theorem for the Cauchy problem (\ref{1.5}) with large initial data for
$\gamma=1+\F{2}{N}$, where $N \geq 5$ odd, with the use of the viscosity method.
The convergence of the Lax-Friedrichs scheme and
the existence of a global solution in $L^\infty$  for large initial data with adiabatic exponent
$ \gamma \in (1, \F{5}{3}]$ were proved  in \cite{8,9}.
In \cite{10}, the global existence of a weak solution was proved for $\gamma \geq 3$ with the
use of the kinetic setting in combination with the compensated compactness method. The method
in \cite{10} was finally improved in \cite{11} to fill the gap $\gamma \in (\F{5}{3},3)$,  and a new proof of the existence of a
global solution for all $\gamma > 1$ was given there. Later on, a new application of the method
in \cite{10} was obtained in \cite{12} on the study of the Euler equations of one-dimensional, compressible
fluid flow, where the linear combinations of weak and strong entropies were invented to replace the weak entropies
used in \cite{7,8,9,10,11}. The isothermal case $ \gamma=1$  with the vacuum was studied in \cite{13}.

Thus, the problem on the existence of a generalized solution of the Cauchy problem (\ref{1.5}) with
bounded initial data (\ref{1.4}) has been completely solved in the case of a polytropic gas.

For the  case of a non-isentropic polytropic gas, namely $ s \ne 0 $,
the  existence theorem of (\ref{1.1}) for small initial data, away from the vacuum,
of locally finite total variation was proved in \cite{14} for $ \gamma \in (1, \F{5}{3}]$  in Lagrangian coordinates, where the proof is based on the finite difference scheme of Glimm \cite{4}.

How to obtain the global existence, for the equations of non-isentropic gas dynamics
 (\ref{1.1}) ( or the simplified system (\ref{1.3}) ) with arbitrarily large
 initial data (\ref{1.4}) including the vacuum, is still a challenging open problem.

 Our aim in this paper is to apply the convergence framework from the compensated compactness theory on the system of isentropic  gas dynamics (\ref{1.5}) to prove the pointwise convergence of the approximation solutions of (\ref{1.3}), and to obtain the global existence of bounded entropy solutions for the Cauchy problem (\ref{1.3}) with the bounded initial data (\ref{1.4}). As a by-product, we obtain a non-classical bounded generalized solution $(\rho,u,s)$,  of the non-isentropic polytropic gas flow (\ref{1.1}), which satisfies the mass equation and the momentum equation, the entropy equation with an extra nonnegative measure in the sense of distributions.

Substituting the first equation in (\ref{1.3}) into the second and the third, we have the following system about the variables $(\rho,u,s)$,
\begin{equation}\left\{
\label{1.6}
\begin{array}{l}
\rho_{t}+u \rho_{x}+ \rho u_{x}=0,         \\\\
       u_t+ u u_x+ \theta^{2} \rho^{\gamma-2} e^{2s} \rho_{x} +  \F{2 \theta^{2}}{ \gamma} \rho^{\gamma-1}  e^{2s} s_{x}=0, \\\\
       s_t+  u s_x=0,
\end{array}\right.
\end{equation}
which, for smooth solutions, is equivalent to system (\ref{1.3}) as well as
system (\ref{1.1}).

Let  the matrix $dF(U)$ of (\ref{1.6}) be
\begin{equation}
 \label{1.7}
dF(U) =\left( \begin{array}{ccc}
u & \rho & 0 \\\\
 \theta^{2} \rho^{\gamma-2} e^{2s} & u & \F{2 \theta^{2}}{ \gamma} \rho^{\gamma-1}  e^{2s} \\\\
0 & 0 & u
\end{array} \right).
\end{equation}
Then three eigenvalues of (\ref{1.3}) are
\begin{equation}
\label{1.8}
\lam_1=u- \theta \rho^{\theta} e^{s} ,  \hspace{0.3cm}
\lam_2=u+ \theta \rho^{\theta} e^{s}, \hspace{0.3cm}
\lam_3= u
\end{equation}
with corresponding right eigenvectors
\begin{equation}
\label{1.9}
r_{1}=(1, -\theta \rho^{\theta-1} e^{s},0)^{T},  \quad
r_{2}=(1, \theta \rho^{\theta-1} e^{s},0)^{T}, \quad
r_{3}=(0, 0,1)^{T}.
\end{equation}
The Riemann invariants of (\ref{1.3}) are functions
$w_{1}( \rho,u,s), w_{2}( \rho,u,s) $ and $w_{3}( \rho,u,s)$ satisfying the equations
\begin{equation}
\label{1.10}
( w_{i \rho}, w_{i u},w_{i s}) \cdot dF= \lam_{i}( w_{i \rho}, w_{i u},w_{i s}), \quad
i=1,2,3.
\end{equation}
Since the system (\ref{1.10}) is not well defined,  we consider $s$ to be a constant, then the Riemann invariants of the isentropic system (\ref{1.5}),
 \begin{equation}
\label{1.11}
w_{1}=z( \rho,u,s)=   \rho^{\theta} e^{s}-u, \hspace{0.3cm}
w_{2}=w( \rho,u,s)= \rho^{\theta} e^{s}+u, \hspace{0.3cm}
\end{equation}
satisfy  the first two equations of (\ref{1.10}).

By simple calculations,
\begin{equation}
\label{1.12}
\nabla \lam_{1} \cdot r_{1}= - \theta(1+\theta) \rho^{\theta-1} e^{s}, \quad
\nabla \lam_{2} \cdot r_{2}=  \theta(1+\theta) \rho^{\theta-1} e^{s}, \quad
\nabla \lam_{3} \cdot r_{3}=0.
\end{equation}
Therefore it follows from (\ref{1.8}) that system (\ref{1.3}) is
strictly hyperbolic in the domain $\{(x,t): \rho(x,t) >0 \} $, while it is hyperbolically degenerate in the domain $ \{(x,t): \rho(x,t) =0 \} $, since
$ \lam_{1}= \lam_{2}=\lam_{3}$ when $ \rho=0$. From (\ref{1.12}), the first two characteristic fields in (\ref{1.3}) are genuinely nonlinear  if the adiabatic exponent $ \gamma \in (1, 3]$, while the system
is no longer genuinely nonlinear at $ \rho=0$ if the adiabatic exponent  $ \gamma >3$; and the third
characteristic field is always linearly degenerate, or of the Temple type \cite{15}.

It is well known that, in order to prove the existence of solutions by using the compensated compactness theory, we should first obtain the a-priori $L^{p} $ estimate, $ 1 < p \leq \infty,$ of the approximate solutions, and
look for enough entropy-entropy flux pairs. Then, we may obtain the measure equations
by applying the div-curl lemma and the representation
of weak limit of solution in terms of Young measure. Finally we must show the
reduction of Young measure to a Dirac measure.

Unfortunately, in the face of nonlinear hyperbolic systems of more than two conservation
laws, we meet the difficulties in all the above three steps.

In fact, except the scalar equation, even if for systems of two
equations, not all pairs of entropy-entropy flux $(\eta,q)$  could
be used to reduce the Young measure to be the Dirac mass since
$\eta_{t}+q_{x}$, where $(\eta,q)$ is a pair of entropy-entropy flux, must be compact in $H^{-1}_{loc}(R \times R^{+})$ when one applies the div-curl lemma
of Tartar \cite{5}.

 When we study the system (\ref{1.3}), the crucial difficulty is how to obtain the a-priori $L^{p}, 1 < p \leq \infty, $ estimates of the approximate solutions because the invariant regions theory \cite{16}, in general, does not work.

 In \cite{17}, the authors studied the following system
\begin{equation}\left\{
\label{1.13}
\begin{array}{l}
\rho_{t}+( \rho u)_{x}=0,         \\\\
       ( \rho(1+s)u)_t+( \rho(1+s)u^2+ P( \rho,s))_x=0, \\\\
       ( \rho s)_t+( \rho u s)_x=0,
       \end{array}\right.
\end{equation}
where $ P( \rho,s)= a^{2} s \F{\rho}{1-\rho}$ and $a$ is a constant. The same as system (\ref{1.3}),  two characteristic fields of (\ref{1.13}) are genuinely nonlinear and the third one is of the Temple type. Under the assumption of a uniform bound on the $L^{\infty} $ norm of the viscosity approximate solutions and other several technical assumptions, the global entropy solution of the Cauchy problem (\ref{1.13}) with bounded initial data was studied with the help of the compactness framework of DiPerna on $ 2 \times 2 $  strictly hyperbolic, genuinely nonlinear systems \cite{18}. Later, the author in \cite{19} showed the existence of invariant regions for the Riemann
problem and obtained global existence using the Glimm scheme \cite{4}.

In \cite{20}, the classical smooth solution of (\ref{1.13}) was obtained when $ P( \rho,s) $ is fixed as $ e^{s} e^{-\F{1}{\rho}}$, where the characteristic fields
are assumed to be nondecreasing.  System (\ref{1.13}) with this special pressure is interesting because it can be diagonalized.

In \cite{21}, the authors studied the following  model of polytropic gas flow with diffusive entropy
\begin{equation}\left\{
\label{1.14}
\begin{array}{l}
\rho_{t}+( \rho u)_{x}=0,         \\\\
       ( \rho u)_t+( \rho u^2+ P( \rho,s))_x=0, \\\\
       (( \rho s)_t+( \rho u s)_x= ( \F{1}{\rho}s_{x})_{x},
\end{array}\right.
\end{equation}
where $P( \rho,s)= e^{(\gamma-1)s} \rho^{\gamma}, \gamma>1.$  With the help of the diffusive term $ ( \F{1}{\rho}s_{x})_{x}$, a skill was used to obtain the a-priori $L^{\infty} $ estimates of the approximate solutions constructed by the  Lax-Friedrichs or Godunov schemes, and the pointwise convergence of the approximate solutions was proved  by using the compactness framework on $ 2 \times 2 $ polytropic gas flow
\cite{8,9,10,11}.

The main contribution of this paper is to obtain the a-priori, $L^{\infty} $
estimates of the viscosity solutions of  (\ref{1.3}).

The classical vanishing viscosity method is to add the viscosity terms
to the right-hand side of system (\ref{1.3}) and consider the Cauchy
problem for the following related parabolic system
\begin{equation}\left\{
\label{1.15}
\begin{array}{l}
\rho_{t}+( \rho u)_{x}= \E \rho_{xx},         \\\\
       ( \rho u)_t+( \rho u^2+ P( \rho) e^{2s})_x= \E (\rho u)_{xx}, \\\\
       ( \rho s)_t+( \rho u s)_x= \E (\rho s)_{xx}
\end{array}\right.
\end{equation}
with bounded initial data
\begin{equation}
\label{1.16} (\rho^{\E},u^{\E},s^{\E})|_{t=0}=(\rho^{\E}_{0}(x),u^{\E}_{0}(x),s^{\E}_{0}(x)), \hspace{0.3cm} \rho^{\E}_{0}(x) \ge \E >0,
\end{equation}
where $ \rho^{\E}_{0}(x)=(\rho_{0}(x)+  \E)*G^{ \E},  u^{\E}_{0}(x)= u_{0}(x)*G^{ \E},  s^{\E}_{0}(x)=s_{0}(x)*G^{ \E}$   are the smooth approximations of
$\rho_{0}(x),u_{0}(x), s_{0}(x)$ and $G^{ \E}$ is a mollifier. However, if we consider
$(\rho,m, \Upsilon)$, where $m= \rho u, \Upsilon= \rho s$ as three
independent variables in (\ref{1.15}), then the
terms $ \rho u^2 = \F{m^2}{\rho}, \rho u s = \F{m}{\rho} \Upsilon $ are singular near the line $ \rho=0$.

Compared with the previous results on (\ref{1.13}) and (\ref{1.14}) introduced above, we mainly need to resolve the following three difficulties when we study the Cauchy problem for System (\ref{1.3}).

{\bf Difficulty I.} How to obtain the positive, lower bound of the
viscosity solutions  $\rho^{\E}$ for the Cauchy problem (\ref{1.15})
and (\ref{1.16})?

To overcome this difficulty, instead of the classical viscosity approximation,
we use again the flux approximation introduced in \cite{22,23} and
consider the following parabolic system
\begin{equation}\left\{
\label{1.17}
\begin{array}{l}
\rho_{t}+( (\rho-2 \delta) u)_{x}= \E \rho_{xx},         \\\\
       ( \rho u)_t+( \rho u^2- \delta u^{2}+ P_{1}(\rho,\delta) e^{2s} )_x= \E (\rho u)_{xx}, \\\\
       ( \rho s)_t+( (\rho-2 \delta) u s)_x= \E (\rho s)_{xx}
\end{array}\right.
\end{equation}
with bounded initial data
\begin{equation}
\label{1.18} (\rho^{\E,\delta},u^{\E,\delta},s^{\E,\delta})|_{t=0}
=(\rho^{\E,\delta}_{0}(x),u^{\E,\delta}_{0}(x),s^{\E,\delta}_{0}(x)),
 \hspace{0.3cm} \rho^{\E,\delta}_{0}(x) \ge 2 \delta >0,
\end{equation}
where $\E, \delta$ are positive small perturbation constants, the
perturbation pressure
\begin{equation}
\label{1.19} P_{1}( \rho, \delta)= \int^{\rho} \F{t-2
\delta}{t}P'(t)dt= k \rho^{\gamma} - 2 \delta k \F{\gamma}{\gamma-1}
\rho^{\gamma-1} ,
\end{equation}
and  $ \rho^{\E,\delta}_{0}(x)=(\rho_{0}(x)+ 2 \delta)*G^{ \E},  u^{\E,\delta}_{0}(x)= u_{0}(x)*G^{ \E},  s^{\E,\delta}_{0}(x)=s_{0}(x)*G^{ \E}$   are the smooth approximations of $\rho_{0}(x),u_{0}(x)$ and $s_{0}(x)$, satisfying $
\lim_{|x| \rightarrow 0} |\F{d^{i}}{dx^{i}}(s^{\E,\delta}_{0}(x))|=0, i=1,2$. Since $ \rho^{\E,\delta}_{0}(x) \geq 2 \delta,$ applying
the maximum principle to the first equation in (\ref{1.17}), we may obtain
 the uniformly positive lower bound $ \rho^{\E,\delta}(x,t) \geq  2 \delta$,
which grantees that $ \rho u^2= \F{m^2}{\rho}, \rho u s = \F{m}{\rho} \Upsilon $ in (\ref{1.17}) are regular. Besides, the flux approximation given in (\ref{1.17}) has the following advantage.

When we consider $s$ as a parameter or $s=0$, (\ref{1.17}) is deduced to the following system of two equations
\begin{equation}\left\{
\label{1.20}
\begin{array}{l}
\rho_{t}+( (\rho-2 \delta) u)_{x}= \E \rho_{xx},         \\\\
       ( \rho u)_t+( \rho u^2- \delta u^{2}+ P_{1}(\rho,\delta)e^{2s}  )_x= \E (\rho u)_{xx}.
\end{array}\right.
\end{equation}
The uniformly lower bound $ \rho^{\E,\delta}(x,t) \geq  2 \delta $  helped us to obtain the proof of the $H^{-1}$ compactness of
$ \eta_{t}+q_{x}$ for any weak entropy-entropy flux pair $(\eta,q)$, and for general pressure function $P(\rho)$ (cf. \cite{22}).

Moreover, system (\ref{1.20}) has the same Riemann invariants and the entropy equation like system (\ref{1.5}).

These special behaviors, of system (\ref{1.17}) as well as (\ref{1.20}), will help us to obtain the uniformly upper bound estimate of $(\rho^{\E, \delta},u^{\E, \delta})$ and to overcome the following difficulty.

{\bf Difficulty II.} How to obtain the uniformly, upper bound $
\rho^{\E, \delta} \leq M$ and $ |u^{\E, \delta}| \leq M$?

The outline to overcome the above difficulty is as follows. First, substituting the first equation in (\ref{1.17}) into the third, we
may rewrite the third equation in (\ref{1.17}) as
\begin{equation}
 \label{1.21}
s_{t}+ \F{(\rho-2 \delta)}{\rho}u s_{x}= \E s_{xx}+ 2 \E \F{\rho_{x}}{\rho} s_{x}.
\end{equation}
By using the technique given in \cite{24,25} (see also
\cite{26}), we can easily prove that $s_{x}^{\E, \delta}$  is uniformly bounded in $L^{1}(R)$.

Second, multiplying $(w_{\rho},w_{m},w_{s}), (z_{\rho},z_{m},z_{s})$ to system (\ref{1.17}), where $ m= \rho u$ and $ w,z$ are given
in (\ref{1.11}), we obtain
\begin{equation}\begin{array}{ll}
\label{1.22}
 w_{t}+ \lam^{\delta}_{2}w_{x}- (\F{2 \delta}{\rho} u + \theta (\rho- 2 \delta) \rho^{ \theta-1} e^{s}) \rho^{ \theta} e^{s}s_{x} \\\\
= \E w_{xx}+ \F{2 \E}{\rho} \rho_{x}w_{x}- \E e^{s} \rho^{\theta-2}( \theta (\theta+1) \rho^{2}_{x}
+2 \theta \rho \rho_{x}s_{x}+  \rho^{2}s^{2}_{x})
\end{array}
\end{equation}
and
\begin{equation}\begin{array}{ll}
\label{1.23}
 z_{t}+ \lam^{\delta}_{1}z_{x}-(\F{2 \delta}{\rho} u - \theta (\rho- 2 \delta) \rho^{ \theta-1} e^{s}) \rho^{ \theta} e^{s}s_{x} \\\\
= \E z_{xx}+ \F{2 \E}{\rho} \rho_{x}z_{x}- \E e^{s} \rho^{\theta-2}( \theta (\theta+1) \rho^{2}_{x}
+2 \theta \rho \rho_{x}s_{x}+  \rho^{2}s^{2}_{x}),
\end{array}
\end{equation}
where
\begin{equation}
\label{1.24}
\lam^{\delta}_{1}= u- \F{\rho-2
\delta}{\rho} \theta \rho^{ \theta} e^{s}, \quad \lam^{\delta}_{2}= u+ \F{\rho-2
\delta}{\rho} \theta \rho^{ \theta} e^{s}
\end{equation}
are first two eigenvalues of the left-hand side of system (\ref{1.17}).

Compared with the case of $s=0$, when we intend to apply the maximum principle to
(\ref{1.22}) and (\ref{1.23}), the functions $ (\F{2 \delta}{\rho} u + \theta (\rho- 2 \delta) \rho^{ \theta-1} e^{s}) \rho^{ \theta} e^{s}s_{x} $
in (\ref{1.22}) and $ (\F{2 \delta}{\rho} u - \theta (\rho- 2 \delta) \rho^{ \theta-1} e^{s}) \rho^{ \theta} e^{s}s_{x}$ in (\ref{1.23}) are two major stumbling blocks.

Fortunately, since $ s_{x} $ is uniformly bounded in $L^{1}(R)$, we might choose a suitable nonnegative, bounded
function $ \beta(x,t), \beta_{x}(x,t) \geq 0,$ and make the transformation of variables
\begin{equation}
\label{1.25}
w=v_{1}+M + \beta(x,t),  \
z=v_{2}+ M - \beta(x,t)
\end{equation}
to obtain
$w_{t}+ \lam^{\delta}_{2}w_{x}=v_{1t}+ \lam^{\delta}_{2}v_{1x}+ (\beta_{t}(x,t)+\lam^{\delta}_{2} \beta_{x}(x,t)),
z_{t}+ \lam^{\delta}_{1}z_{x}=v_{2t}+ \lam^{\delta}_{1}v_{2x}- (\beta_{t}(x,t)+\lam^{\delta}_{1} \beta_{x}(x,t))$, where the extra functions $ \beta_{t}(x,t)+\lam^{\delta}_{i} \beta_{x}(x,t), i=1,2,$ could be used to control $(\F{2 \delta}{\rho} u + \theta (\rho- 2 \delta) \rho^{ \theta-1} e^{s}) \rho^{ \theta} e^{s}s_{x} $ and $ (\F{2 \delta}{\rho} u - \theta (\rho- 2 \delta) \rho^{ \theta-1} e^{s}) \rho^{ \theta} e^{s}s_{x}$.

If we could obtain the estimates $ v_{1} \leq 0, v_{2} \leq 0$ by using the maximum
principle to the new variables $(v_{1}, v_{2})$, clearly from the transformation (\ref{1.25}), the bounds  $w(\rho^{\E, \delta}, u^{\E, \delta}) \leq M + \beta(x,t),
z(\rho^{\E, \delta}, u^{\E, \delta}) \leq M - \beta(x,t) $  follow on.

{\bf Difficulty III.} How to prove the pointwise convergence
of $(\rho^{\E, \delta}, u^{\E, \delta},s^{\E, \delta}) $ as $\E, \delta$ go to zero?

First, as introduced in \cite{24},  the  pointwise convergence of $
s^{\E, \delta}$ can be obtained by using the div-curl lemma to some special pairs of functions $(c, F(s))$, where $c$ is a constant and $F(s)$ is a suitable function of $s$ since
the $L^{1}(R)$ estimate of $ s_{x}^{\E, \delta}(\cdot, t)$ and the  compactness of $c_{t}+ F(s^{\E,
\delta})_{x}$ in $H^{-1}_{loc} $.

Second, to prove the pointwise convergence
of $(\rho^{\E, \delta}, u^{\E, \delta}) $ as $\E, \delta$ go to zero, we may fix the variable $s$ or
think of $s$ as a parameter,  with the help of the compactness framework \cite{8,9,10,11} on the $ 2 \times 2 $ polytropic gas dynamic system (\ref{1.20}),
we may prove the pointwise convergence of $(\rho^{\E, \delta}, u^{\E, \delta}) $ and obtain the global existence of solutions.

The main results of this paper are listed in the following Theorem 1.1.
\begin{theorem}
Let the initial data $ (\rho_{0}(x),u_{0}(x),s_{0}(x)) $ be bounded in $L^{\infty}(\mathbb{R});
|s_{0}(x)| \leq N, |s_{0x}|_{L^{1}(\mathbb{R})} \leq c_{0} < 1 $, for two positive constants
$N, c_{0}$. Then, (I) for fixed $ \E, \delta$, the global smooth solution, $(\rho^{\E, \delta}, u^{\E, \delta}, s^{\E, \delta})$ of the Cauchy problem
(\ref{1.17}) and (\ref{1.18}), exists and satisfies
\begin{equation}
\label{1.26}
|s^{\E, \delta}| \leq N, \ |s_{x}^{\E, \delta}(\cdot,t)| _{L^{1}(\mathbb{R})} \leq c_{0} <1,
\end{equation}
and
\begin{equation}\left\{
\label{1.27}
\begin{array}{l}
 z(\rho^{\E, \delta},u^{\E, \delta}, s^{\E, \delta}) 
 \leq M- c \int_{-\infty}^{x} |s_{x}^{\E, \delta}| dx \leq M,
\\\\
 w(\rho^{\E, \delta},u^{\E, \delta}, s^{\E, \delta}) \leq M+ c \int_{-\infty}^{x} |s_{x}^{\E, \delta}| dx
 \leq M+ c c_{0},
\end{array}\right.
\end{equation}
where $z,w$ are the Riemann invariants, of isentropic system  (\ref{1.5}),  given in (\ref{1.11}) and  $c, M$, are two suitable large constants,
satisfying  $M \leq c $ and $ 0 < c c_{0} < M $.

(II) There exists a subsequence of
$(\rho^{\E, \delta}(x,t),u^{\E, \delta}(x,t),s^{\E, \delta}(x,t)),$ which converges pointwisely,  on the set $\rho_{+}= \{(x,t): \rho(x,t) >0 \}$, to a set of bounded functions $(\rho(x,t),u(x,t), s(x,t)) $ as $
\E, \delta $ tend to zero, and the limit is a weak entropy solution
of the Cauchy problem (\ref{1.3})-(\ref{1.4}).
\end{theorem}
\begin{define} A set of bounded functions $(\rho(x,t),u(x,t),s(x,t))$  is called a weak entropy
solution of the Cauchy problem (\ref{1.3})-(\ref{1.4}) if
\begin{equation}\left\{
\label{1.28}
\begin{array}{l}
      \int_{0}^{\infty} \int_{- \infty}^{\infty} \rho \phi_{t}+ \rho u \phi_{x}  \phi dxdt
      + \int_{- \infty}^{\infty} \rho_{0}(x) \phi(x,0) dx=0,  \\  \\
       \int_{0}^{\infty} \int_{- \infty}^{\infty}  \rho u \phi_t+( \rho u^2+ \F{ \theta^{2}}{\gamma} \rho^{\gamma} e^{2s}) \phi_x
      dxdt
      + \int_{- \infty}^{\infty} \rho_{0}(x) u_{0}(x) \phi(x,0) dx=0,  \\\\
      \int_{0}^{\infty} \int_{- \infty}^{\infty} \rho s \phi_{t}+ \rho u s \phi_{x}  dxdt
      + \int_{- \infty}^{\infty} \rho_{0}(x) s_{0}(x) \phi(x,0) dx=0,
\end{array}\right.
\end{equation}
holds for all test function $ \phi \in C_{0}^{1}(R \times R^{+})$
and
\begin{equation}
\label{1.29}
\int_0^{ \infty} \int_{ -\infty}^{\infty} \eta(\rho,m,\Upsilon)
\phi_t+ q( \rho,m,\Upsilon) \phi_x  dxdt \geq 0
\end{equation}
holds for any nonnegative test function $ \phi \in C_0^{ \infty}(R
\times R^{+}- \{t=0\}),$ where  $ m= \rho u, \Upsilon= \rho s $ and
$ ( \eta, q), \eta(0,m,\Upsilon)=0, $ is a
pair of convex, weak entropy-entropy flux of system (\ref{1.3}).
\end{define}
In the next section, we will introduce a technique from the maximum principle,
to set up a relation between the Riemann invariants of isentropic system and the entropy variable $s$,
to study the uniform estimates of the approximate viscosity solutions of the parabolic system (\ref{1.17}) with the initial data (\ref{1.18}). Under the conditions in  Theorem 1.1, we may  obtain the  estimates (\ref{1.26}), (\ref{1.27})  on  $(w^{\E, \delta}(x,t) $ and $z^{\E, \delta}(x,t))$. Then, based on these
estimates, we obtain the pointwise convergence of $(\rho^{\E, \delta}(x,t), u^{\E, \delta}(x,t), s^{\E, \delta}(x,t))$  by applying  the compensated compactness
 theory  on the polytropic gas dynamics (\ref{1.5})  (the first two equations in (\ref{1.3})).

The study of (\ref{1.3}) could be considered as the beginning
of a study of the non-isentropic polytropic gas flow (\ref{1.1}).

In fact, beside the standard viscosity terms, if we add the extra perturbation function
$A(x,t)$ to
(\ref{1.1}) and consider the following parabolic system
\begin{equation}\left\{
\label{1.30}
\begin{array}{l}
\rho_{t}+( \rho u)_{x}= \E \rho_{xx},         \\\\
       ( \rho u)_{t}+( \rho u^2+ k \rho^{\gamma} e^{2s})_x= \E ( \rho u)_{xx}, \\\\
       (\F{1}{2} \rho u^2+ \F{k}{\gamma-1} \rho^{\gamma} e^{2s} )_{t}+
       (u( \F{1}{2} \rho u^2+ \F{k \gamma}{\gamma-1} \rho^{\gamma} e^{2s}))_{x}
       \\\\=\E (\F{1}{2} \rho u^2+ \F{k}{\gamma-1} \rho^{\gamma} e^{2s} )_{xx}-
       \E A(x,t),
\end{array}\right.
\end{equation}
where
\begin{equation}
\label{1.31}
A(x,t) =  \rho u^{2}_{x}
+ k  \rho^{\gamma-2} e^{2s}( \gamma \rho^{2}_{x} +
\F{4 }{\gamma-1} \rho^{2} s^{2}_{x} + 4 \rho \rho_{x} s_{x}),
\end{equation}
then we may prove that (\ref{1.15}) and (\ref{1.30}) are completely same.

To prove the equivalent of (\ref{1.15}) and (\ref{1.30}), we
substitute the first equation in (\ref{1.30}) into the second to obtain
\begin{equation}\begin{array}{ll}
\label{1.32}
u_{t}+ u u_{x}+ k \gamma \rho^{\gamma-2} e^{2s}  \rho_{x}
+2 k  \rho^{\gamma-1} e^{2s} s_{x} = \E u_{xx}+ 2 \E \F{\rho_{x}}{\rho} u_{x}.
\end{array}
\end{equation}
By simple calculations,
\begin{equation}\begin{array}{ll}
\label{1.33}
(\F{1}{2} \rho u^{2}+ \F{k}{\gamma-1} \rho^{\gamma} e^{2s} )_{t}+
       (u( \F{1}{2} \rho u^2+ \F{k \gamma}{\gamma-1} \rho^{\gamma} e^{2s})_{x} \\\\
       = \F{1}{2}u^{2} \rho_{t}+  + \rho u u_{t} + \F{k \gamma}{\gamma-1} \rho^{\gamma-1} e^{2s} \rho_{t} + \F{2 k}{\gamma-1} \rho^{\gamma} e^{2s}s_{t} \\\\+ \F{1}{2}u^{2} (\rho u)_{x}+ \rho u (u u_{x})+
        \rho u (  k \gamma \rho^{\gamma-2} e^{2s}  \rho_{x}) \\\\
        +  \F{k \gamma}{\gamma-1} \rho^{\gamma-1} e^{2s} (\rho u)_{x}
        + \F{2 k \gamma}{\gamma-1} \rho^{\gamma} e^{2s}u s_{x} \\\\
        = \F{1}{2}u^{2} (\rho_{t}+ (\rho u)_{x} )
        + \rho u (u_{t}+ u u_{x}+ k \gamma \rho^{\gamma-2} e^{2s}  \rho_{x}
        + 2 k  \rho^{\gamma-1} e^{2s} s_{x}) \\\\
        + \F{k \gamma}{\gamma-1} \rho^{\gamma-1} e^{2s} (\rho_{t}+ (\rho u)_{x} )
        + \F{2 k}{\gamma-1} \rho^{\gamma} e^{2s}( s_{t} + u s_{x}) =B(x,t)
\end{array}
\end{equation}
and
\begin{equation}\begin{array}{ll}
\label{1.34}
\E (\F{1}{2} \rho u^2+ \F{k}{\gamma-1} \rho^{\gamma} e^{2s} )_{xx}
= \F{1}{2}  u^{2} \E \rho_{xx}+ \E (2 u \rho_{x} u_{x}+ \rho u u_{xx}) + \E \rho u^{2}_{x} \\\\ + \F{k \gamma}{\gamma-1} \rho^{\gamma-1} e^{2s} \E \rho_{xx}
+ \F{2 k}{\gamma-1} \rho^{\gamma} e^{2s} \E s_{xx}
+ \E k \gamma \rho^{\gamma-2} e^{2s} \rho^{2}_{x} \\\\
+  \E \F{4 k \gamma}{\gamma-1} \rho^{\gamma-1} e^{2s} \rho_{x} s_{x}
+ \E \F{4 k }{\gamma-1} \rho^{\gamma} e^{2s} s^{2}_{x} =C(x,t).
\end{array}
\end{equation}
From the third equation in (\ref{1.30}), we have $ B(x,t)=C(x,t)- \E A(x,t)$,
and thus the following equation
\begin{equation}\begin{array}{ll}
\label{1.35}
 \F{2 k}{\gamma-1} \rho^{\gamma} e^{2s}( s_{t} + u s_{x})
= \E \rho u^{2}_{x}
+ \F{2 k}{\gamma-1} \rho^{\gamma} e^{2s} \E s_{xx}
+ \E k \gamma \rho^{\gamma-2} e^{2s} \rho^{2}_{x} \\\\
+  \E \F{4 k \gamma}{\gamma-1} \rho^{\gamma-1} e^{2s} \rho_{x} s_{x}
+ \E \F{4 k }{\gamma-1} \rho^{\gamma} e^{2s} s^{2}_{x}- \E A(x,t),
\end{array}
\end{equation}
where we used the first equation in (\ref{1.30})
and (\ref{1.32}).

From the first equation in (\ref{1.30}), we have
\begin{equation}
\label{1.36}
 \F{2 k}{\gamma-1} \rho^{\gamma-1} e^{2s} (s \rho_{t}+s( \rho u)_{x})
 = \E \rho_{xx} \F{2 k}{\gamma-1} \rho^{\gamma-1} e^{2s}s.
\end{equation}
Summing up (\ref{1.35}) and (\ref{1.36}), we have
\begin{equation}\begin{array}{ll}
\label{1.37}
 \F{2 k}{\gamma-1} \rho^{\gamma-1} e^{2s}( (\rho s)_{t} + (\rho u s)_{x})
= \F{2 k}{\gamma-1} \rho^{\gamma-1} e^{2s} \E (\rho s)_{xx} \\\\
+ \E \rho u^{2}_{x}
+ \E k  \rho^{\gamma-2} e^{2s}( \gamma \rho^{2}_{x} +
\F{4 }{\gamma-1} \rho^{2} s^{2}_{x} + 4 \rho \rho_{x} s_{x}) - \E A(x,t)
 \\\\
 = \F{2 k}{\gamma-1} \rho^{\gamma-1} e^{2s} \E (\rho s)_{xx}
\end{array}
\end{equation}
which deduces the third equation in (\ref{1.15}).

Similarly, we may prove the following parabolic system
\begin{equation}\left\{
\label{1.38}
\begin{array}{l}
\rho_{t}+( (\rho-2 \delta) u)_{x}= \E \rho_{xx},         \\\\
       ( \rho u)_t+( \rho u^2- \delta u^{2}+ P_{1}(\rho,\delta) e^{2s} )_x= \E (\rho u)_{xx}, \\\\
        (\F{1}{2} \rho u^2+ \F{k}{\gamma-1} \rho^{\gamma} e^{2s} )_{t}+
       (u( \F{1}{2} \rho u^2+ \F{k \gamma}{\gamma-1} \rho^{\gamma} e^{2s}))_{x}
       \\\\=\E (\F{1}{2} \rho u^2+ \F{k}{\gamma-1} \rho^{\gamma} e^{2s} )_{xx}-
       \E A(x,t)+ \delta D(x,t),
\end{array}\right.
\end{equation}
 is equivalent to
the flux-viscosity approximate system (\ref{1.17}), where
\begin{equation}
\label{1.39}
D(x,t) = \F{1}{3} (u^{3})_{x} +  \F{4 k}{\gamma-1} \rho^{\gamma-1} e^{2s} u s_{x} + 2 k \F{\gamma}{\gamma-1} (\rho^{\gamma-1} e^{2s}u)_x.
\end{equation}
In fact, we rewrite the first and the second equations in (\ref{1.38}) as
follows
\begin{equation}
\label{1.40}
\rho_{t}+(\rho u)_{x}-2 \delta u_{x}= \E \rho_{xx}
\end{equation}
and
\begin{equation}
\label{1.41}
( \rho u)_t+( \rho u^2+ k \rho^{\gamma} e^{2s} )_x
- 2 \delta u u_{x}- 2 \delta k \F{\gamma}{\gamma-1} (\rho^{\gamma-1} e^{2s} )_x
= \E (\rho u)_{xx}.
\end{equation}
Accordingly, (\ref{1.32}) is rewritten as follows
\begin{equation}\begin{array}{ll}
\label{1.42}
u_{t}+ u u_{x}+ k \gamma \rho^{\gamma-2} e^{2s}  \rho_{x}
+2 k  \rho^{\gamma-1} e^{2s} s_{x}- 2 \delta k \F{\gamma}{(\gamma-1) \rho} (\rho^{\gamma-1} e^{2s} )_x  \\\\
 = \E u_{xx}+ 2 \E \F{\rho_{x}}{\rho} u_{x}.
\end{array}
\end{equation}
By using the calculations in (\ref{1.33})-(\ref{1.34}), the equations
(\ref{1.40})-(\ref{1.42}), and the third equation in the system (\ref{1.38}),
we obtain
\begin{equation}\begin{array}{ll}
\label{1.43}
\F{1}{2} u^{2} ( 2 \delta u_{x})+ \rho u ( 2 \delta k \F{\gamma}{(\gamma-1) \rho} (\rho^{\gamma-1} e^{2s} )_x ) + k \F{\gamma}{(\gamma-1)} \rho^{\gamma-1} e^{2s}
( 2 \delta u_{x}) \\\\+
 \F{2 k}{\gamma-1} \rho^{\gamma} e^{2s}( s_{t} + u s_{x})
= \E \rho u^{2}_{x}
+ \F{2 k}{\gamma-1} \rho^{\gamma} e^{2s} \E s_{xx}
+ \E k \gamma \rho^{\gamma-2} e^{2s} \rho^{2}_{x} \\\\
+  \E \F{4 k \gamma}{\gamma-1} \rho^{\gamma-1} e^{2s} \rho_{x} s_{x}
+ \E \F{4 k }{\gamma-1} \rho^{\gamma} e^{2s} s^{2}_{x}- \E A(x,t) + \delta D(x,t).
\end{array}
\end{equation}
Multiplying $ \F{2 k}{\gamma-1} \rho^{\gamma-1} e^{2s} s $ to the first
equation in (\ref{1.38}), we have
\begin{equation}
\label{1.44}
 \F{2 k}{\gamma-1} \rho^{\gamma-1} e^{2s} (s \rho_{t}+s( (\rho- 2 \delta) u)_{x})
 = \E \rho_{xx} \F{2 k}{\gamma-1} \rho^{\gamma-1} e^{2s}s.
\end{equation}
Since the following terms in (\ref{1.43})
\begin{equation}\begin{array}{ll}
\label{1.45}
 \rho u ( 2 \delta k \F{\gamma}{(\gamma-1) \rho} (\rho^{\gamma-1} e^{2s} )_x ) + k \F{\gamma}{(\gamma-1)} \rho^{\gamma-1} e^{2s}
( 2 \delta u_{x}) \\\\
= 2 \delta k \F{\gamma}{\gamma-1} (\rho^{\gamma-1} e^{2s}u)_x,
\end{array}
\end{equation}
\begin{equation}\begin{array}{ll}
\label{1.46}
\F{2 k}{\gamma-1} \rho^{\gamma} e^{2s}( s_{t} + u s_{x})  \\\\
= \F{2 k}{\gamma-1} \rho^{\gamma-1} e^{2s}( \rho s_{t} +(\rho- 2 \delta) u s_{x})
+ \delta \F{4 k}{\gamma-1} \rho^{\gamma-1} e^{2s} u s_{x}
\end{array}
\end{equation}
and
\begin{equation}
\label{1.47}
 \F{1}{2} u^{2} ( 2 \delta u_{x})= \F{\delta}{3} (u^{3})_{x},
\end{equation}
we may obtain the following equation by summing up  (\ref{1.43}) and (\ref{1.44})
\begin{equation}\begin{array}{ll}
\label{1.48}
\F{\delta}{3} (u^{3})_{x} + \delta \F{4 k}{\gamma-1} \rho^{\gamma-1} e^{2s} u s_{x} + 2 \delta k \F{\gamma}{\gamma-1} (\rho^{\gamma-1} e^{2s}u)_x \\\\
+
 \F{2 k}{\gamma-1} \rho^{\gamma-1} e^{2s}( (\rho s)_{t} + ((\rho- 2 \delta) u s)_{x})
= \F{2 k}{\gamma-1} \rho^{\gamma-1} e^{2s} \E (\rho s)_{xx} + \delta D(x,t)
\end{array}
\end{equation}
which deduces the third equation in (\ref{1.17}).

{\bf Remark 1.} From the analysis above,  under the conclusions given in
Theorem 1.1, since $|s_{x}|$ is bounded in $L^{1}_{loc}(R \times R^{+})$,
we may prove from (\ref{1.38}) that
\begin{equation}
\label{1.49}
   \E A(x.t)= \E \rho u^{2}_{x}
+ \E k  \rho^{\gamma-2} e^{2s}( \gamma \rho^{2}_{x} +
\F{4 }{\gamma-1} \rho^{2} s^{2}_{x} + 4 \rho \rho_{x} s_{x})
\end{equation}
is bounded in $ L^{1}_{loc}(R \times R^{+}) $, and hence converges weakly
to a nonnegative integrable measure $ \mu(x,t) \geq 0$, as $ \E, \delta \rightarrow 0$, and the limit $(\rho,u,s)$,
given in Theorem 1.1, satisfies
\begin{equation}\left\{
\label{1.50}
\begin{array}{l}
      \int_{0}^{\infty} \int_{- \infty}^{\infty} \rho \phi_{t}+ \rho u \phi_{x}  \phi dxdt
      + \int_{- \infty}^{\infty} \rho_{0}(x) \phi(x,0) dx=0,  \\  \\
       \int_{0}^{\infty} \int_{- \infty}^{\infty}  \rho u \phi_t+( \rho u^2+ \F{ \theta^{2}}{\gamma} \rho^{\gamma} e^{2s}) \phi_x
      dxdt
      + \int_{- \infty}^{\infty} \rho_{0}(x) u_{0}(x) \phi(x,0) dx=0,
  \end{array}\right.
\end{equation}
for all test function $ \phi \in C_{0}^{1}(R \times R^{+})$, and
 \begin{equation}\begin{array}{ll}
\label{1.51}
      \int_{0}^{\infty} \int_{- \infty}^{\infty} (\F{1}{2} \rho u^2+ \F{k}{\gamma-1} \rho^{\gamma} e^{2s}) \phi_{t}+ (u( \F{1}{2} \rho u^2+ \F{k \gamma}{\gamma-1} \rho^{\gamma} e^{2s})) \phi_{x} dxdt
      \\\\
      + \int_{- \infty}^{\infty} (\F{1}{2} \rho_{0}(x) (u_{0}(x))^2+ \F{k}{\gamma-1} (\rho_{0}(x))^{\gamma} e^{2s_{0}(x)}) \phi(x,0) dx\\\\
      = \int_{0}^{\infty} \int_{- \infty}^{\infty}  \mu(x,t)  \phi dxdt \geq 0,
\end{array}
\end{equation}
for any nonnegative test function $ \phi \in C_{0}^{1}(R \times R^{+})$.

This means that, if we could prove $ \mu (x,t)=0 $, then the global existence
of the classical entropy solutions, in the sense of Definition 1, for the system of non-isentropic polytropic gas flow (\ref{1.1}) could be obtained by using the vanishing viscosity method introduced in (\ref{1.15}) or (\ref{1.17}).
Clearly, for classical smooth solutions, Systems (\ref{1.1}) and (\ref{1.3})
are equivalent because $ \mu (x,t)=0$.

{\bf Remark 2.} The inequality (\ref{1.51}) is natural because it is the entropy
inequality, for the physical entropy-entropy flux pair $(\eta_{0},q_{0})$
of the isentropic gas dynamics (\ref{1.5}), in the sense of (\ref{1.29}), where $s$ is considered as a parameter, and
\begin{equation}
\label{1.52}
(\eta_{0},q_{0})=(\F{1}{2} \rho u^2+ \F{k}{\gamma-1} \rho^{\gamma} e^{2s}, u( \F{1}{2} \rho u^2+ \F{k \gamma}{\gamma-1} \rho^{\gamma} e^{2s})).
\end{equation}
So, we obtain the following definition of a non-classical generalized solution.

\begin{define} A set of bounded functions $(\rho(x,t),u(x,t),s(x,t))$  is called a non-classical generalized solution of the Cauchy problem (\ref{1.1}) and (\ref{1.4}) if the pair of functions
$(\rho(x,t),u(x,t))$ satisfies the mass and the momentum equations in (\ref{1.1}) in the sense of distributions, and  the function $ s(x,t) $ satisfies the entropy equation in (\ref{1.1}) with a nonnegative integrable measure $\mu(x,t)$ in the sense of (\ref{1.51}).
\end{define}
Based on the conclusions in Theorem 1.1, we have the following existence
of solutions for the Cauchy problem (\ref{1.1}) and (\ref{1.4}).
\begin{theorem}
Under the conditions of the initial data in Theorem 1.1, the
Cauchy problem of the non-isentropic polytropic gas flow (\ref{1.1}) and
the initial data (\ref{1.4}) has a non-classical generalized solution.
\end{theorem}

\section{ \bf Proof of Theorem 1.1.}

In this section we shall prove Theorem 1.1. First, following the standard theory of semilinear parabolic systems, the local existence result of the Cauchy problem (\ref{1.17}), (\ref{1.18}) can be easily
obtained by applying the contraction mapping principle to an integral representation for a solution. Second, we have the following lemma about the estimates on $ s^{\E, \delta}(x, t)$.
\begin{lemma}
If $  |s_{0}(x)| \leq N$ and $  |s_{0x}|_{L^{1}(\mathbb{R})} \leq c_{0}$, then the functions $ s^{\E,\delta}(x,t)$ obtained from the Cauchy problem
(\ref{1.17}) and (\ref{1.18}) satisfy
\begin{equation}
 \label{2.1}
 | s^{\E,\delta}(x,t)| \leq N,  \ |s^{\E,\delta}_{x}(\cdot,t)|_{L^{1}(\mathbb{R})} \leq c_{0}.
\end{equation}
Moreover,
\begin{equation}
 \label{2.2}
\int_{- \infty}^{x}
|\omega|_{t} dx+ \F{(\rho-2 \delta)}{\rho}u |\omega|= \E |\omega|_{x}+
2 \E \F{\rho_{x}}{\rho} |\omega|,
\end{equation}
where $ \omega(x,t)=s_{x}^{\E, \delta}(x, t)$ and $N, c_{0}<1$ are two positive constants.
\end{lemma}
{ \bf Proof of Lemma 2.1.} The first assertion in (\ref{2.1}) can be obtained by applying the maximum principle to (\ref{1.21}) directly.

Differentiating Equation (\ref{1.21}) with respect to $x$, we have
\begin{equation}
 \label{2.3}
\omega_{t}+ (\F{(\rho-2 \delta)}{\rho}u)_{x} \omega+
\F{(\rho-2 \delta)}{\rho}u \omega_{x}= \E \omega_{xx}+
2 \E (\F{\rho_{x}}{\rho})_{x} \omega + 2 \E \F{\rho_{x}}{\rho} \omega_{x}.
\end{equation}
The  second assertion in (\ref{2.1}) can be proved by the methods given in
\cite{1,24,26}.

For any fixed $t$, since $s^{\E, \delta}(x, t)$ is of bounded total variation,
$s_{x}^{\E, \delta}(x, t)$ changes the sign at most the countable points $x_{i}, i=1,2,3...$. Thus the following equality is true except at these points $( x_{i},t)$
  \begin{equation}
 \label{2.4}
|\omega|_{t}+ (\F{(\rho-2 \delta)}{\rho}u |\omega|)_{x} = \E |\omega|_{xx}+
2 \E (\F{\rho_{x}}{\rho}|\omega|)_{x}.
\end{equation}
Integrating (\ref{2.4}) in $ (-\infty,x)$,  we get the proof of (\ref{2.2}). {\bf Lemma 2.1 is proved}.

To obtain the a-priori upper estimates of $ w^{\E, \delta}(x,t)$ and $z^{\E, \delta}(x,t))$ given in (\ref{1.27}), we first prove (\ref{1.22}) and (\ref{1.23}).

{ \bf The proof of (\ref{1.22}) and (\ref{1.23}).} By simple calculations,
\[ \left\{
\begin{array}{l}
w_{\rho}= \theta \rho^{\theta-1}e^{s}- \F{m}{\rho^{2}} , \hspace{0.2cm}  w_{m}=\F{1}{\rho}, \hspace{0.2cm}  w_{s}=\rho^{\theta}e^{s},  \\\\
w_{\rho \rho}= \theta (\theta-1) \rho^{\theta-2}e^{s}+2 \F{m}{\rho^{3}}, \hspace{0.2cm}  w_{\rho m}=- \F{1}{\rho^{2}}, \hspace{0.2cm}  w_{\rho s}=\theta \rho^{\theta-1}e^{s},  \\\\
w_{mm}=0, \hspace{0.2cm} w_{sm}=0, \hspace{0.2cm} w_{ss}=\rho^{\theta}e^{s},  \\\\
      z_{\rho}= \theta \rho^{\theta-1}e^{s}+ \F{m}{\rho^{2}} , \hspace{0.2cm}  z_{m}=-\F{1}{\rho}, \hspace{0.2cm}  z_{s}=\rho^{\theta}e^{s},  \\\\
z_{\rho \rho}= \theta (\theta-1) \rho^{\theta-2}e^{s}-2 \F{m}{\rho^{3}}, \hspace{0.2cm}  z_{\rho m}= \F{1}{\rho^{2}}, \hspace{0.2cm}  z_{\rho s}=\theta \rho^{\theta-1}e^{s},  \\\\
z_{mm}=0, \hspace{0.2cm} z_{sm}=0, \hspace{0.2cm} z_{ss}=\rho^{\theta}e^{s}.
      \end{array}\right. \]
Multiplying $(w_{\rho},w_{m})$  to the first two equations in system (\ref{1.17}), and $ w_{s}$ to
(\ref{1.21}), then adding the results,  we obtain an equality whose left-hand side is
\[ \begin{array}{ll}
L= w_{t}+ \lam^{\delta}_{2}w_{x} - \lam^{\delta}_{2}w_{s}s_{x} + \F{ \rho- 2 \delta}{ \rho} u w_{s}s_{x} \\\\
= w_{t}+ \lam^{\delta}_{2}w_{x} - (\F{2 \delta}{\rho} u + \theta (\rho- 2 \delta) \rho^{ \theta-1} e^{s}) \rho^{ \theta} e^{s}s_{x}
\end{array}  \]
and the right-hand side is
\[ \begin{array}{ll}
R= \E w_{xx}- \E (w_{\rho \rho} \rho_{x}^{2}+ 2 w_{\rho m} \rho_{x}m_{x}+ w_{m m} m_{x}^{2}
+ w_{s s} s_{x}^{2}+  2 w_{s m} s_{x}m_{x}\\\\
+ 2 w_{\rho s } \rho_{x}s_{x})
+ 2 \E \F{ \rho_{x} }{ \rho} w_{s}s_{x}  \\\\
=\E w_{xx}+ \F{2 \E}{\rho} \rho_{x}w_{x}- \E e^{s} \rho^{\theta-2}( \theta (\theta+1) \rho^{2}_{x}
+2 \theta \rho \rho_{x}s_{x}+  \rho^{2}s^{2}_{x}),
\end{array} \]
so we proved (\ref{1.22}). Similarly,  if we multiply $(z_{\rho},z_{m})$  to the first two equations in system (\ref{1.17}), $ z_{s}$ to (\ref{1.21}) and add the results, we may obtain the proof of (\ref{1.23}).
\begin{lemma}
Make the transformation of variables
 \begin{equation}
\label{2.5}
 w=v_{1}+M + c \int_{-\infty}^{x} |s_{x}| dx,  \
z=v_{2}+ M - c \int_{-\infty}^{x} |s_{x}| dx,
\end{equation}
 where $M, c $ are suitable large constants, satisfying $ 0 < M \leq c, \
  c c_{0} < M $.  Then the new variables $v_{1}, v_{2}$ satisfy the following
system of two inequalities
\begin{equation}\left\{
\label{2.6}
\begin{array}{l}
      v_{1t}+a_{1}(x,t) v_{1x}+ b_{1}(x,t) v_{1}+c_{1}(x,t)v_{2}  \leq \E v_{1xx},    \\  \\
       v_{2t}+a_{2}(x,t)v_{2x}+b_{2}(x,t)v_{2}+c_{2}(x,t)v_{1}  \leq \E v_{2xx},
\end{array}\right.
\end{equation}
where $ a_{i}(x,t),b_{i}(x,t), c_{i}(x,t) \leq 0, i=1,2, $ are suitable functions.
\end{lemma}
{ \bf Proof of Lemma 2.2.} Using the transformation (\ref{2.5}), we have from (\ref{1.22}) that
\begin{equation}\begin{array}{ll}
\label{2.7}
 v_{1t}+c \int_{-\infty}^{x} |s_{x}|_{t} dx+ \lam^{\delta}_{2}v_{1x}+ \lam^{\delta}_{2} c |s_{x}| - (\F{2 \delta}{\rho} u + \theta (\rho- 2 \delta) \rho^{ \theta-1} e^{s}) \rho^{ \theta} e^{s}s_{x} \\\\
= \E v_{1xx} + \E c|s_{x}|_{x} + \F{2 \E}{\rho} \rho_{x}v_{1x}+ \F{2 \E}{\rho} \rho_{x} c |s_{x}|
- \E e^{s} \rho^{\theta-2}( \theta (\theta+1) \rho^{2}_{x}
+2 \theta \rho \rho_{x}s_{x}+  \rho^{2}s^{2}_{x}) \\\\
\leq \E v_{1xx} + \E c|s_{x}|_{x} + \F{2 \E}{\rho} \rho_{x}v_{1x}+ \F{2 \E}{\rho} \rho_{x} c |s_{x}|.
\end{array}
\end{equation}
From (\ref{2.2}), we have
\begin{equation}
 \label{2.8}
c \int_{-\infty}^{x} |s_{x}|_{t} dx- \E c |s_{x}|_{x}- 2 c \E \F{\rho_{x}}{\rho} |s_{x}| =- c \F{(\rho-2 \delta)}{\rho}u |s_{x}|.
\end{equation}
Since, we finally obtain the estimates $v_{1} \leq 0, v_{2} \leq 0,$ which deduce
from (\ref{2.5}) that
 \begin{equation}\begin{array}{ll}
\label{2.9}
w \leq M + c \int_{-\infty}^{x} |s_{x}| dx,  \ z \leq M - c \int_{-\infty}^{x} |s_{x}|dx,
\   \rho^{ \theta} e^{s}  = \F{1}{2} (w+z) \leq  M,
\end{array}
\end{equation}
thus it is enough if we may prove Lemma 2.2 in the region $ \rho^{ \theta} e^{s} \leq M$.

By simple calculations, 
\begin{equation}\begin{array}{ll}
\label{2.10}
\lam^{\delta}_{2} c |s_{x}| - (\F{2 \delta}{\rho} u + \theta (\rho- 2 \delta) \rho^{ \theta-1} e^{s}) \rho^{ \theta} e^{s}s_{x}- c \F{(\rho-2 \delta)}{\rho}u |s_{x}|\\\\
= (u+ \F{\rho-2
\delta}{\rho} \theta \rho^{ \theta} e^{s}) c|s_{x}| - (\F{2 \delta}{\rho} u + \theta (\rho- 2 \delta) \rho^{ \theta-1} e^{s}) \rho^{ \theta} e^{s}s_{x}- c \F{(\rho-2 \delta)}{\rho}u |s_{x}|,
\end{array}
\end{equation}
which we write as $I_{1}(x,t)$. First, at the points $(x,t)$, where $s_{x} \geq 0$, we have
$ s_{x}= |s_{x}| $ and
\begin{equation}\begin{array}{ll}
\label{2.11}
I_{1}(x,t) = (u+ \F{\rho-2
\delta}{\rho} \theta \rho^{ \theta} e^{s}) c|s_{x}| - (\F{2 \delta}{\rho} u + \theta (\rho- 2 \delta) \rho^{ \theta-1} e^{s}) \rho^{ \theta} e^{s}|s_{x}|- c \F{(\rho-2 \delta)}{\rho}u |s_{x}|
\\\\
= \theta \F{\rho-2 \delta}{\rho} \rho^{ \theta} e^{s} |s_{x}| (c-\rho^{ \theta} e^{s})
+ \F{2 \delta}{\rho} |s_{x}| u (c-\rho^{ \theta} e^{s})
\geq \F{2 \delta}{\rho} |s_{x}| u (c-\rho^{ \theta} e^{s}) \\\\
= \F{2 \delta}{\rho} |s_{x}|  (c-\rho^{ \theta} e^{s}) \F{1}{2}(v_{1}-v_{2}+ 2 c \int_{-\infty}^{x} |s_{x}| dx)
\geq \F{2 \delta}{\rho} |s_{x}|  (c-\rho^{ \theta} e^{s}) \F{1}{2}(v_{1}-v_{2}),
\end{array}
\end{equation}
where we used $ \rho^{ \theta} e^{s} \leq M \leq c$. Similarly, at the
points $(x,t)$, where $s_{x} \leq 0$, we have $ s_{x}= -|s_{x}| $ and
\begin{equation}\begin{array}{ll}
\label{2.12}
I_{1}(x,t) = (u+ \F{\rho-2
\delta}{\rho} \theta \rho^{ \theta} e^{s}) c|s_{x}| + (\F{2 \delta}{\rho} u + \theta (\rho- 2 \delta) \rho^{ \theta-1} e^{s}) \rho^{ \theta} e^{s}|s_{x}|- c \F{(\rho-2 \delta)}{\rho}u |s_{x}|
\\\\
= \theta \F{\rho-2 \delta}{\rho} \rho^{ \theta} e^{s} |s_{x}| (c+\rho^{ \theta} e^{s})
+ \F{2 \delta}{\rho} |s_{x}| u (c+\rho^{ \theta} e^{s})
\geq \F{2 \delta}{\rho} |s_{x}| u (c+\rho^{ \theta} e^{s}) \\\\
= \F{2 \delta}{\rho} |s_{x}|  (c+\rho^{ \theta} e^{s}) \F{1}{2}(v_{1}-v_{2}+ 2 c \int_{-\infty}^{x} |s_{x}| dx)
\geq \F{2 \delta}{\rho} |s_{x}|  (c+\rho^{ \theta} e^{s}) \F{1}{2}(v_{1}-v_{2}),
\end{array}
\end{equation}
Thus, we have from (\ref{2.7}), (\ref{2.8}), (\ref{2.11}) and (\ref{2.12})  that
\begin{equation}\begin{array}{ll}
\label{2.13}
 v_{1t}+ ( \lam^{\delta}_{2}- \F{2 \E}{\rho} \rho_{x}  ) v_{1x}+
 \F{2 \delta}{\rho} |s_{x}|  (c- sgn (s_{x}) \rho^{ \theta} e^{s}) \F{1}{2}(v_{1}-v_{2})
 \leq  \E v_{1xx},
\end{array}
\end{equation}
which gives us the proof of the first inequality in (\ref{2.6}).

Similarly,  we have from (\ref{1.23}) and the transformation (\ref{2.5}) that
\begin{equation}\begin{array}{ll}
\label{2.14}
 v_{2t}-c \int_{-\infty}^{x} |s_{x}|_{t} dx+ \lam^{\delta}_{1}v_{2x}- \lam^{\delta}_{1} c|s_{x}| - (\F{2 \delta}{\rho} u - \theta (\rho- 2 \delta) \rho^{ \theta-1} e^{s}) \rho^{ \theta} e^{s}s_{x} \\\\
= \E v_{2xx} - \E c|s_{x}|_{x}+ \F{2 \E}{\rho} \rho_{x}v_{2x}- \F{2 \E}{\rho} \rho_{x} c|s_{x}|
- \E e^{s} \rho^{\theta-2}( \theta (\theta+1) \rho^{2}_{x}
+2 \theta \rho \rho_{x}s_{x}+  \rho^{2}s^{2}_{x}) \\\\
\leq \E v_{2xx} - \E c|s_{x}|_{x}+ \F{2 \E}{\rho} \rho_{x}v_{2x}- \F{2 \E}{\rho} \rho_{x} c|s_{x}|.
\end{array}
\end{equation}
By simple calculations, 
\begin{equation}\begin{array}{ll}
\label{2.15}
-\lam^{\delta}_{1} c |s_{x}| - (\F{2 \delta}{\rho} u - \theta (\rho- 2 \delta) \rho^{ \theta-1} e^{s}) \rho^{ \theta} e^{s}s_{x}+ c \F{(\rho-2 \delta)}{\rho}u |s_{x}|\\\\
= (-u+ \F{\rho-2
\delta}{\rho} \theta \rho^{ \theta} e^{s}) c|s_{x}| - (\F{2 \delta}{\rho} u - \theta (\rho- 2 \delta) \rho^{ \theta-1} e^{s}) \rho^{ \theta} e^{s}s_{x}+ c \F{(\rho-2 \delta)}{\rho}u |s_{x}|=I_{2}(x,t).
\end{array}
\end{equation}
At the points  $(x,t)$, where $s_{x} \geq 0$, we have
$ s_{x}= |s_{x}| $ and
\begin{equation}\begin{array}{ll}
\label{2.16}
I_{2}(x,t) = \theta \F{\rho-2 \delta}{\rho} \rho^{ \theta} e^{s} |s_{x}| (c+\rho^{ \theta} e^{s})
- \F{2 \delta}{\rho} |s_{x}| u (c+\rho^{ \theta} e^{s})
= I(x,t) |s_{x}| (c+\rho^{ \theta} e^{s}),
\end{array}
\end{equation}
where
\begin{equation}\begin{array}{ll}
\label{2.17}
I(x,t)= \theta \F{\rho-2 \delta}{\rho} \rho^{ \theta} e^{s} - \F{2 \delta}{\rho} u
= \theta \F{\rho-2 \delta}{\rho} \rho^{ \theta} e^{s}
+ \F{2 \delta}{\rho} (z- \rho^{ \theta} e^{s})  \\\\
= \theta \F{\rho-2 \delta}{\rho} \rho^{ \theta} e^{s}
+ \F{2 \delta}{\rho}(v_{2}+M-c \int_{-\infty}^{x} |s_{x}| dx - \rho^{ \theta} e^{s}) \\\\
= (\theta-(\theta+1)\F{2 \delta}{\rho})  \rho^{ \theta} e^{s}
+ \F{2 \delta}{\rho} (M-c \int_{-\infty}^{x} |s_{x}| dx) + \F{2 \delta}{\rho} v_{2}.
\end{array}
\end{equation}
Similarly, at the points  $(x,t)$, where $s_{x} \leq 0$, we have
$ s_{x}= -|s_{x}| $ and
\begin{equation}\begin{array}{ll}
\label{2.18}
I_{2}(x,t) = \theta \F{\rho-2 \delta}{\rho} \rho^{ \theta} e^{s} |s_{x}| (c-\rho^{ \theta} e^{s})
- \F{2 \delta}{\rho} |s_{x}| u (c-\rho^{ \theta} e^{s})
= I(x,t) |s_{x}| (c-\rho^{ \theta} e^{s}).
\end{array}
\end{equation}
Now, we analyze the function $I(x,t)$. First, at the points $(x,t)$, where $\theta-(\theta+1)\F{2 \delta}{\rho(x,t)} \geq 0$,
we have $ I(x,t) \geq \F{2 \delta}{\rho} v_{2} $ due to
$ c|s_{x}|_{L^{1}} \leq c c_{0} < M$.

Second, at the points $(x,t)$, where $\theta-(\theta+1)\F{2 \delta}{\rho(x,t)} \leq 0$
or $ \rho(x,t) \leq 2 \delta \F{\theta+1}{\theta} $, we have
\begin{equation}\begin{array}{ll}
\label{2.19}
(\theta-(\theta+1)\F{2 \delta}{\rho})  \rho^{ \theta} e^{s}
+ \F{2 \delta}{\rho} (M-c \int_{-\infty}^{x} |s_{x}| dx)  \\\\
\geq (\theta-(\theta+1)\F{2 \delta}{\rho}) ( 2 \delta \F{\theta+1}{\theta} )^{ \theta} e^{s}+ \F{2 \delta}{\rho} (M-c \int_{-\infty}^{x} |s_{x}| dx) \\\\
= ( \theta \F{\rho-2 \delta}{\rho} - \F{2 \delta}{\rho}) ( 2 \delta \F{\theta+1}{\theta} )^{ \theta} e^{s} + \F{2 \delta}{\rho} (M-c \int_{-\infty}^{x} |s_{x}| dx)  \\\\
\geq - \F{2 \delta}{\rho} ( 2 \delta \F{\theta+1}{\theta} )^{ \theta} e^{s} + \F{2 \delta}{\rho} (M-c \int_{-\infty}^{x} |s_{x}| dx) \\\\
= \F{2 \delta}{\rho} (M-c \int_{-\infty}^{x} |s_{x}| dx- ( 2 \delta \F{\theta+1}{\theta} )^{ \theta} e^{s} )
\geq 0
\end{array}
\end{equation}
if we choose $ \delta $ to be sufficiently small. Therefore we obtain the second inequality in (\ref{2.6}) from (\ref{2.14}), (\ref{2.2}), (\ref{2.15})-(\ref{2.19}),
and complete the proof of Lemma 2.2.

Since the initial data are bounded, then at $ t=0, v_{1}(x,0) = w-M-c  \int_{-\infty}^{x} |s_{x}| dx \leq 0,
v_{2}(x,0) = w-M+c \int_{-\infty}^{x} |s_{x}| dx \leq 0 $ for large $M$. Applying the maximum principle to (\ref{2.6}), we have $ v_{1}(x,t) \leq 0, v_{2}(x,t) \leq 0 $ for any time $t >0$. Thus we have $w \leq M+ c \int_{-\infty}^{x} |s_{x}| dx,
z \leq M- c \int_{-\infty}^{x} |s_{x}| dx $, which deduce the estimates in
(\ref{1.27}) and so the estimates on the viscosity solutions
$(\rho^{\E,\delta}(x,t),u^{\E,\delta}(x,t),s^{\E,\delta}(x,t)) $ of
the Cauchy problem (\ref{1.17}) and (\ref{1.18}):
\begin{equation}\begin{array}{ll}
\label{2.20}
2 \delta \leq \rho^{\E,\delta}(x,t) \leq N, \ |u^{\E,\delta}(x,t)| \leq N, \
 |s^{\E,\delta}(x,t)| \leq N, \ |s_{x}^{\E,\delta}(\cdot,t)|_{L^{1}} \leq c_{0} <1,
\end{array}
\end{equation}
where $N$ is a positive constant depending only the bound of the initial data,
but being independent of $\E$ and $ \delta$.

With the uniformly bounded estimates in (\ref{2.20}), we may extend the local
solution of the Cauchy problem (\ref{1.17}) and (\ref{1.18}) step by step, until
an arbitrary large time $T$ and obtain the global existence of solution. So, the part {\bf (I)} in Theorem 1.1 is proved.

To obtain the pointwise convergence of a subsequence of $(\rho^{\E,\delta},u^{\E,\delta},s^{\E,\delta})$, we first have the following lemma
\begin{lemma}
There exists a subsequence (still labelled
 $s^{\E,\delta}(x,t))$ such that
\begin{equation}
\label{2.21}
 s^{\E,\delta}(x,t) \rightarrow s(x,t)
\end{equation}
almost everywhere on the set $\rho_{+}= \{(x,t): \rho(x,t) >0 \}$, where
$\rho(x,t)$ is the weak-star limit of $ \rho^{\E,\delta}(x,t)$.
\end{lemma}
{\bf Proof of Lemma 2.3.}
Since $ s^{\E,\delta}_{x} $  and $ ((s^{\E,\delta})^{2})_{x} $ are
uniformly bounded in $W_{loc}^{-1, \infty}(R \times R^{+}) \cap L^{1}_{loc}(R \times R^{+})$, then
\begin{equation}
 \label{2.22}
 c_{t}+ s^{\E,\delta}_{x}, \ c_{t}+ ((s^{\E,\delta})^{2})_{x} \quad
 \mbox{ are compact in } \quad  H^{-1}_{loc}(R \times R^{+})
\end{equation}
by Murat's Lemma \cite{6}, where $c$ is a constant.

 Multiplying $(\F{ \partial \eta_{0}}{ \partial \rho},\F{ \partial \eta_{0}}{ \partial m})$ to the first two equations in System (\ref{1.17}), $ \F{ \partial \eta_{0}}{ \partial s}$ to (\ref{1.21}), where
 $ \eta_{0}$ is given in (\ref{1.52}), then adding the results, we have
 \begin{equation}\begin{array}{ll}
\label{2.23}
        (\F{1}{2} \rho u^2+ \F{k}{\gamma-1} \rho^{\gamma} e^{2s} )_{t}+
       (u( \F{1}{2} \rho u^2+ \F{k \gamma}{\gamma-1} \rho^{\gamma} e^{2s}))_{x}
       =\E (\F{1}{2} \rho u^2+ \F{k}{\gamma-1} \rho^{\gamma} e^{2s} )_{xx}
       \\\\- \E (\rho u^{2}_{x} + k  \rho^{\gamma-2} e^{2s}( \gamma \rho^{2}_{x} +
\F{4 }{\gamma-1} \rho^{2} s^{2}_{x} + 4 \rho \rho_{x} s_{x})) \\\\
+ \delta (\F{1}{3} (u^{3})_{x} +  \F{4 k}{\gamma-1} \rho^{\gamma-1} e^{2s} u s_{x} + 2 k \F{\gamma}{\gamma-1} (\rho^{\gamma-1} e^{2s}u)_x).
\end{array}
\end{equation}
Let $K \subset R \times R^{+} $ be an arbitrary compact set and choose $ \phi \in C_{0}^{\infty}(R \times R^+) $ such that $ \phi_{K}=1, 0 \leq \phi \leq 1 $.

Multiplying Equation (\ref{2.23}) by $ \phi$ and integrating over $R \times R^+ $, we may obtain
\begin{equation}\begin{array}{ll}
\label{2.24}
\int_{0}^{\infty} \int_{-\infty}^{\infty} \E (\rho u^{2}_{x} + k  \rho^{\gamma-2} e^{2s}( \gamma \rho^{2}_{x} +
\F{4 }{\gamma-1} \rho^{2} s^{2}_{x} + 4 \rho \rho_{x} s_{x}))  \phi dxdt
\leq M( \phi),
\end{array}
\end{equation}
where we used the $L^{1}$ local integrability of $S_{x}$, and hence that
\begin{equation}
\label{2.25}
\E \rho u_{x}^2, \ \E \rho^{\gamma} s_{x}^2,
\ \E \rho^{\gamma-2} \rho_{x}^2 \mbox{ are bounded in } L^{1}_{loc}(R \times R^+),
\end{equation}
due to
\begin{equation}
\label{2.26}
\gamma \rho^{2}_{x} + \F{4 }{\gamma-1} \rho^{2} s^{2}_{x} + 4 \rho \rho_{x} s_{x} \geq
c (\rho^{2}_{x} + \rho^{2} s^{2}_{x})
\end{equation}
for a suitable constant $ c>0$.

For any $\varphi\in H^{1}_{0}(R \times R^{+})$, we have from the estimates in (\ref{2.25}) that
\begin{equation}\begin{array}{ll}
\label{2.27}
|\int_{0}^{\infty} \int_{- \infty}^{\infty} \E \rho_{xx} \varphi dx dt|
= |\int_{0}^{\infty} \int_{- \infty}^{\infty} \E \rho_{x} \varphi_{x} dx dt| \\\\
\leq (\int_{0}^{\infty} \int_{- \infty}^{\infty} \E \rho^{\gamma-2} \rho^{2}_{x} |\varphi_{x}| dx dt)^{\F{1}{2}} (\int_{0}^{\infty} \int_{- \infty}^{\infty} \E \rho^{2-\gamma} |\varphi_{x}| dx dt)^{\F{1}{2}} \\\\
\leq M  (\int_{0}^{\infty} \int_{- \infty}^{\infty} \E \rho^{2-\gamma} |\varphi_{x}| dx dt)^{\F{1}{2}}
\rightarrow 0
\end{array}
\end{equation}
and
\begin{equation}\begin{array}{ll}
\label{2.28}
|\int_{0}^{\infty} \int_{- \infty}^{\infty} \E (\rho s)_{xx} \varphi dx dt|
= |\int_{0}^{\infty} \int_{- \infty}^{\infty} \E (\rho_{x} s+ \rho s_{x}) \varphi_{x} dx dt| \\\\
\leq (\int_{0}^{\infty} \int_{- \infty}^{\infty} \E \rho^{\gamma-2} \rho^{2}_{x} |\varphi_{x}| dx dt)^{\F{1}{2}} (\int_{0}^{\infty} \int_{- \infty}^{\infty} \E \rho^{2-\gamma} s^{2} |\varphi_{x}| dx dt)^{\F{1}{2}} \\\\+
(\int_{0}^{\infty} \int_{- \infty}^{\infty} \E \rho^{\gamma} s^{2}_{x} |\varphi_{x}| dx dt)^{\F{1}{2}} (\int_{0}^{\infty} \int_{- \infty}^{\infty} \E \rho^{2-\gamma} |\varphi_{x}| dx dt)^{\F{1}{2}} \\\\
\leq M  (\int_{0}^{\infty} \int_{- \infty}^{\infty} \E \rho^{2-\gamma} |\varphi_{x}| dx dt)^{\F{1}{2}}
\rightarrow 0,
\end{array}
\end{equation}
because we may choose $\E$ to go zero much faster than $\delta$ such that
$ \E \rho^{2-\gamma} \rightarrow 0 $ as $\E, \delta$ go to zero. Then
we have from the first and the third equations in (\ref{1.17}) that
\begin{equation}\begin{array}{ll}
\label{2.29}
\rho^{\E,\delta}_{t}+((\rho^{\E,\delta}-2 \delta) u^{\E,\delta})_x \ \mbox{ and } \ ( \rho^{\E,\delta} s^{\E,\delta})_t+( (\rho^{\E,\delta}-2 \delta) u^{\E,\delta} s^{\E,\delta})_x
\end{array}
\end{equation}
are compact in $H^{-1}_{loc}(R \times R^{+})$.

Thus we may apply the div-curl lemma to the pairs of functions
\begin{equation}
 \label{2.30}
(c, s^{\E,\delta}), \quad (\rho^{\E,\delta},(\rho^{\E,\delta}-2 \delta) u^{\E,\delta})
\end{equation}
and
\begin{equation}
 \label{2.31}
(c, s^{\E,\delta}), \quad ( \rho^{\E,\delta} s^{\E,\delta},(\rho^{\E,\delta}-2 \delta) u^{\E,\delta} s^{\E,\delta})
\end{equation}
respectively to obtain
\begin{equation}
 \label{2.32}
\overline{\rho^{\E,\delta}} \cdot \overline{s^{\E,\delta}}=
\overline{ \rho^{\E,\delta} s^{\E,\delta}}, \ \mbox{ and } \
\overline{s^{\E,\delta}} \cdot \overline{\rho^{\E,\delta} s^{\E,\delta}}=
\overline{ \rho^{\E,\delta} (s^{\E,\delta})^{2}},
\end{equation}
where $ \overline{f(\theta^{\E,\delta})}$ denotes the weak-star
limit of $f(\theta^{\E,\delta})$.

Let $ (\overline{\rho^{\E,\delta}},\overline{s^{\E,\delta}})=(\rho,s)$. We have
from (\ref{2.32}) that
\begin{equation}
 \label{2.33}
\overline{ \rho^{\E,\delta} (s^{\E,\delta}-s)^{2}}=\overline{ \rho^{\E,\delta}
(s^{\E,\delta})^{2}}-2s \overline{ \rho^{\E,\delta}s^{\E,\delta}}+ \rho
s^{2}=0.
\end{equation}
Furthermore, we may apply the div-curl lemma to the pair of functions
\begin{equation}
 \label{2.34}
(c, (s^{\E,\delta})^{2}), \quad (\rho^{\E,\delta},(\rho^{\E,\delta}-2 \delta) u^{\E,\delta})
\end{equation}
to obtain
\begin{equation}
 \label{2.35}
\overline{\rho^{\E,\delta}} \cdot \overline{ (s^{\E,\delta})^{2}}=
\overline{ \rho^{\E,\delta} (s^{\E,\delta})^{2}}.
\end{equation}
Using (\ref{2.32}),(\ref{2.33}) and (\ref{2.35}), we have
\begin{equation}\begin{array}{ll}
 \label{2.36}
\overline{ \rho (s^{\E,\delta}-s)^{2}}=
\rho \overline{  (s^{\E,\delta})^{2}}- 2 \rho s \overline{s^{\E,\delta}}
+ \rho s^{2} \\\\=
\overline{ \rho^{\E,\delta}
(s^{\E,\delta})^{2}}-2s \overline{ \rho^{\E,\delta}s^{\E,\delta}}+ \rho
s^{2}=0.
\end{array}
\end{equation}
Then
\begin{equation}
 \label{2.37}
  \rho (s^{\E,\delta}-s)^{2} \rightarrow 0, \quad \mbox{ a.e.},
  \end{equation}
which deduces the proof of Lemma 2.3.

After we have the pointwise convergence of $ s^{\E,\delta} $, we may consider
$s$ as a constant (or a parameter), and study the following system
\begin{equation}\left\{
\label{2.38}
\begin{array}{l}
\rho_{t}+( (\rho-2 \delta) u)_{x}=0,         \\\\
       ( \rho u)_t+( \rho u^2- \delta u^{2}+ P_{1}(\rho,\delta)e^{2s}  )_x=0.
\end{array}\right.
\end{equation}
For smooth solutions, system (\ref{2.38}) is equivalent to the following system
\begin{equation}\left\{
\label{2.39}
\begin{array}{l}
      \rho_{t}+(-2 \delta u+ \rho u)_{x}=0         \\
       u_t+( \F{1}{2}u^2+ \displaystyle \int_{2 \delta}^\rho
\F{(t-2 \delta)P'(t)}{t^{2}} dt e^{2 s})_x=0,
\end{array}\right.
\end{equation}
and particularly, both systems have the same entropy-entropy flux pairs.
Thus any entropy-entropy flux pair
$( \eta( \rho, u,s), q( \rho,u,s))$ of
system (\ref{2.39}) satisfies the additional system
\begin{equation}
\label{2.40}
q_{ \rho}= u \eta_{\rho} + \F{(\rho-2 \delta) P'( \rho) e^{2s}}{ \rho^{2}} \eta_u,
\quad
 q_{u}=( \rho-2 \delta) \eta_{\rho} +u \eta_u.
\end{equation}
Eliminating the $q$ from (\ref{2.40}), we have
\begin{equation}
\label{2.41}
\eta_{ \rho \rho}=   \theta^{2} \rho^{\gamma-3} e^{2s} \eta_{uu}.
\end{equation}
Therefore, system (\ref{2.38}) has the same entropy equation, as system (\ref{1.5}), given in \cite{6}.

An entropy $ \eta(\rho,u,s) $ of system (\ref{2.38})
is called a weak entropy if $\eta(0,u,s)=0$, that is,
a solution of Equation (\ref{2.41}) with the
special initial conditions:
\begin{equation}
\label{2.42}
 \eta( \rho=0,u,s)=0,
\quad
 \eta_{ \rho}( \rho=0,u,s)=f(u,0)=g(u),
\end{equation}
where $g(u)$ is an arbitrary given function of $u$.
The solution of (\ref{2.41})-(\ref{2.42}) is
given by the following lemma:
\begin{lemma} For $ \rho \geq 0, u, w \in R,$
let
\begin{equation}
\label{2.43}
G( \rho,s,w)=( \rho^{\gamma-1} e^{2s}-w^2)_+^{ \lam},\quad
\lam= \F{3- \gamma}{2( \gamma-1)},
\end{equation}
where the notation $x_+= \sup (0,x)$. Then we have
\begin{equation}\begin{array}{ll}
\label{2.44}
\eta( \rho,u,s)&= \displaystyle \int_{R} g( \xi) G( \rho,s, \xi-u) d \xi \\\\
&= \rho \displaystyle \int_0^1[ \tau(1- \tau)]^{ \lam}
g(u+ \rho^{ \theta}e^{s}- 2 \rho^{ \theta}e^{s} \tau) d \tau;
\end{array}
\end{equation}
and the weak entropy flux $q(\rho,u,s)$ of system (\ref{2.38}) associated with $ \eta(\rho,u,s)$ is
\begin{equation}\begin{array}{ll}
\label{2.45}
q(\rho,u,s) = \rho \int_0^1[ \tau(1- \tau)]^{ \lam}
g(u+ \rho^{ \theta} e^{s}- 2 \rho^{ \theta}e^{s} \tau)(u+ \theta(1-2 \tau)
 \rho^{ \theta}e^{s}) d \tau \\\\
 - 2 \delta  \int_0^1[ \tau(1- \tau)]^{ \lam}
G(u+ \rho^{ \theta}e^{s}- 2 \rho^{ \theta}e^{s} \tau) d \tau  \\\\
- 2 \delta  \theta
 \int_0^1[ \tau(1- \tau)]^{ \lam}
g(u+ \rho^{ \theta} e^{s}- 2 \rho^{ \theta}e^{s} \tau)(1-2 \tau)
 \rho^{ \theta}e^{s} d \tau,
\end{array}
\end{equation}
where $ G(y)= \int^{y} g(x) dx$.
 \end{lemma}
 {\bf Proof of Lemma 2.4.}
 The weak entropy formula (\ref{2.44}) is given in \cite{10} (See also \cite{27} or Lemma 8.2.1 in
\cite{28}).

Using the second equation in (\ref{2.40}) and the
weak solution formula (\ref{2.44}), we have
\begin{equation}\begin{array}{ll}
\label{2.46}
q_{u}(\rho,u,s) = \eta +  \theta
 \int_0^1[ \tau(1- \tau)]^{ \lam}
g'(u+ \rho^{ \theta} e^{s}- 2 \rho^{ \theta} e^{s} \tau)(1-2 \tau)
 \rho^{ \theta+1}e^{s} d \tau  \\\\
+  u \rho  \int_0^1[ \tau(1- \tau)]^{ \lam}
g'(u+ \rho^{ \theta} e^{s}- 2 \rho^{ \theta}e^{s} \tau) d \tau \\\\
- 2 \delta  \int_0^1[ \tau(1- \tau)]^{ \lam}
g(u+ \rho^{ \theta}e^{s}- 2 \rho^{ \theta}e^{s} \tau) d \tau  \\\\
- 2 \delta  \theta
 \int_0^1[ \tau(1- \tau)]^{ \lam}
g'(u+ \rho^{ \theta} e^{s}- 2 \rho^{ \theta}e^{s} \tau)(1-2 \tau)
 \rho^{ \theta}e^{s} d \tau.
\end{array}
\end{equation}
Since
\begin{equation}
\begin{array}{ll}
\label{2.47}
   \int^u ug'(u+ \rho^{ \theta}e^{s}- 2 \rho^{ \theta}e^{s} \tau)du
\\\\ = ug(u+ \rho^{ \theta}e^{s}- 2 \rho^{ \theta}e^{s} \tau)-
 \int^u g(u+ \rho^{ \theta}e^{s}- 2 \rho^{ \theta}e^{s} \tau)du,
\end{array}
\end{equation}
we get from (\ref{2.46}) that
\begin{equation}\begin{array}{ll}
\label{2.48}
q(\rho,u,s) =  u \eta
 + \theta
 \int_0^1[ \tau(1- \tau)]^{ \lam}
g(u+ \rho^{ \theta} e^{s}- 2 \rho^{ \theta}e^{s} \tau)(1-2 \tau)
 \rho^{ \theta+1}e^{s} d \tau \\\\
 - 2 \delta  \int_0^1[ \tau(1- \tau)]^{ \lam}
G(u+ \rho^{ \theta}e^{s}- 2 \rho^{ \theta}e^{s} \tau) d \tau  \\\\
- 2 \delta  \theta
 \int_0^1[ \tau(1- \tau)]^{ \lam}
g(u+ \rho^{ \theta} e^{s}- 2 \rho^{ \theta} e^{s} \tau)(1-2 \tau)
 \rho^{ \theta}e^{s} d \tau \\\\
 = q_{1}(\rho,u,s) - 2 \delta q_{2}(\rho,u,s),
\end{array}
\end{equation}
where $ G(y)= \int^{y} g(x) dx$,
 \begin{equation}
 \begin{array}{ll}
\label{2.49}
q_{1}(\rho,u,s)= u \eta
 + \theta
 \int_0^1[ \tau(1- \tau)]^{ \lam}
g(u+ \rho^{ \theta} e^{s}- 2 \rho^{ \theta} e^{s} \tau)(1-2 \tau)
 \rho^{ \theta+1}e^{s} d \tau \\\\
 =\rho \int_0^1[ \tau(1- \tau)]^{ \lam}
g(u+ \rho^{ \theta} e^{s}- 2 \rho^{ \theta} e^{s} \tau)(u+ \theta(1-2 \tau)
 \rho^{ \theta}e^{s}) d \tau
\end{array}
\end{equation}
and
\begin{equation}\begin{array}{ll}
\label{2.50}
q_{2}(\rho,u,s)= \int_0^1[ \tau(1- \tau)]^{ \lam}
G(u+ \rho^{ \theta}e^{s}- 2 \rho^{ \theta}e^{s} \tau) d \tau  \\\\
+  \theta
 \int_0^1[ \tau(1- \tau)]^{ \lam}
g(u+ \rho^{ \theta} e^{s}- 2 \rho^{ \theta}e^{s} \tau)(1-2 \tau)
 \rho^{ \theta}e^{s} d \tau.
\end{array}
\end{equation}
To complete the proof of Lemma 6, we still need to prove that $ q(\rho,u,s) $, given in (\ref{2.48}), satisfies the first equation in (\ref{2.40}), namely
\begin{equation}
\label{2.51}
q_{1 \rho}= u \eta_{\rho} + \theta^{2} \rho^{\gamma-2} e^{2s} \eta_u,
\quad
 - 2 \delta q_{2 \rho}=- 2 \delta  \theta^{2} \rho^{\gamma-3} e^{2s} \eta_u.
\end{equation}
By simple calculations,
\begin{equation}
 \begin{array}{ll}
\label{2.52}
q_{1 \rho}(\rho,u,s)= u \eta_{\rho} +
\theta^{2}
 \int_0^1[ \tau(1- \tau)]^{ \lam}
g'(u+ \rho^{ \theta} e^{s}- 2 \rho^{ \theta}e^{s} \tau)(1-2 \tau)^{2}
 \rho^{2 \theta}e^{2s} d \tau  \\\\
+ \theta (\theta+1)
 \int_0^1[ \tau(1- \tau)]^{ \lam}
g(u+ \rho^{ \theta} e^{s}- 2 \rho^{ \theta}e^{s} \tau)(1-2 \tau)
 \rho^{ \theta}e^{s} d \tau  \\\\
 = u \eta_{\rho} + \theta^{2} \rho^{2 \theta} e^{2s}
 \int_0^1[ \tau(1- \tau)]^{ \lam}
g'(u+ \rho^{ \theta} e^{s}- 2 \rho^{ \theta}e^{s} \tau)  d \tau  \\\\
 + \theta^{2}
 \int_0^1[ \tau(1- \tau)]^{ \lam}
g'(u+ \rho^{ \theta} e^{s}- 2 \rho^{ \theta} e^{s} \tau) 4 \tau (\tau -1)
 \rho^{2 \theta}e^{2s} d \tau  \\\\
+ \theta (\theta+1)
 \int_0^1[ \tau(1- \tau)]^{ \lam}
g(u+ \rho^{ \theta} e^{s}- 2 \rho^{ \theta}e^{s} \tau)(1-2 \tau)
 \rho^{ \theta}e^{s} d \tau \\\\
 = u \eta_{\rho} + \theta^{2} \rho^{\gamma-2} e^{2s} \eta_u
 + 2 \theta^{2} \rho^{ \theta}e^{s}
 \int_0^1[ \tau(1- \tau)]^{ \lam+1}
 d (g(u+ \rho^{ \theta} e^{s}- 2 \rho^{ \theta}e^{s} \tau)) \\\\
 + \theta (\theta+1)
 \int_0^1[ \tau(1- \tau)]^{ \lam}
g(u+ \rho^{ \theta} e^{s}- 2 \rho^{ \theta}e^{s} \tau)(1-2 \tau)
 \rho^{ \theta}e^{s} d \tau \\\\
 =  u \eta_{\rho} + \theta^{2} \rho^{\gamma-2} e^{2s} \eta_u
\end{array}
\end{equation}
due to $ \theta(\theta+1)= 2 \theta (\lam+1)$. Moreover,
\begin{equation}
 \begin{array}{ll}
\label{2.53}
q_{2 \rho}(\rho,u,s)=\theta \int_0^1[ \tau(1- \tau)]^{ \lam}
g(u+ \rho^{ \theta}e^{s}- 2 \rho^{ \theta}e^{s} \tau) (1-2 \tau)
 \rho^{ \theta-1}e^{s} d \tau  \\\\
 +   \theta^{2} \int_0^1[ \tau(1- \tau)]^{ \lam}
g'(u+ \rho^{ \theta} e^{s}- 2 \rho^{ \theta}e^{s} \tau)(1-2 \tau)^{2}
 \rho^{2 \theta-1}e^{2 s} d \tau \\\\
+  \theta^{2} \int_0^1[ \tau(1- \tau)]^{ \lam}
g(u+ \rho^{ \theta} e^{s}- 2 \rho^{ \theta}e^{s} \tau)(1-2 \tau)
 \rho^{ \theta-1}e^{s} d \tau \\\\
 =(\theta+\theta^{2}) \int_0^1[ \tau(1- \tau)]^{ \lam}
g(u+ \rho^{ \theta}e^{s}- 2 \rho^{ \theta}e^{s} \tau) (1-2 \tau)
 \rho^{ \theta-1}e^{s} d \tau  \\\\
 +   \theta^{2} \int_0^1[ \tau(1- \tau)]^{ \lam}
g'(u+ \rho^{ \theta} e^{s}- 2 \rho^{ \theta}e^{s} \tau)(1-2 \tau)^{2}
 \rho^{2 \theta-1}e^{2 s} d \tau.
\end{array}
\end{equation}
Since
\begin{equation}
 \begin{array}{ll}
\label{2.54}
\theta^{2}  \int_0^1[ \tau(1- \tau)]^{ \lam}
g'(u+ \rho^{ \theta} e^{s}- 2 \rho^{ \theta}e^{s} \tau)(1-2 \tau)^{2}
 \rho^{2 \theta-1}e^{2 s} d \tau \\\\
 = \theta^{2} \rho^{2 \theta-1}e^{2 s} \int_0^1[ \tau(1- \tau)]^{ \lam}
g'(u+ \rho^{ \theta} e^{s}- 2 \rho^{ \theta}e^{s} \tau) d \tau \\\\
 - 4 \theta^{2}  \int_0^1[ \tau(1- \tau)]^{ \lam+1}
g'(u+ \rho^{ \theta} e^{s}- 2 \rho^{ \theta}e^{s} \tau)
 \rho^{2 \theta-1}e^{2 s} d \tau \\\\
 = \theta^{2} \rho^{\gamma-3} e^{2s} \eta_u - 2 \theta^{2} (\lam+1)
 \rho^{ \theta-1}e^{s} \int_0^1[ \tau(1- \tau)]^{ \lam}
g(u+ \rho^{ \theta} e^{s}- 2 \rho^{ \theta}e^{s} \tau) d \tau,
\end{array}
\end{equation}
we have the second equality in
(\ref{2.51}) by summing up (\ref{2.53}) and (\ref{2.54}),
\begin{equation}
\label{2.55}
 q_{2 \rho}=\theta^{2} \rho^{\gamma-3} e^{2s} \eta_u,
\end{equation}
and hence obtain the proof of Lemma 2.4.
\begin{lemma}
\begin{equation}
\label{2.56}
 \eta_{t} (\rho^{\E,\delta},u^{\E,\delta},s^{\E,\delta})
 +q_{1x}(\rho^{\E,\delta},u^{\E,\delta},s^{\E,\delta})
 \quad
 \mbox{ are compact in} \quad  H^{-1}_{loc}(R \times
R^{+}),
\end{equation}
with respect to the viscosity solutions $(\rho^{\E,\delta},u^{\E,\delta},s^{\E,\delta})$ of the Cauchy problem
(\ref{1.17}) and (\ref{1.18}), where $ \eta, q_{1}$ are given in (\ref{2.44}) and (\ref{2.49}).
\end{lemma}
{\bf Proof of Lemma 2.5.} Let $ \Phi(\rho,u,s,\tau)= u+ \rho^{ \theta}e^{s}- 2 \rho^{ \theta}e^{s} \tau$. By simple calculations, we have from (\ref{2.44}) that
\begin{equation}\left\{
\label{2.57}
\begin{array}{l}
   \eta_{\rho}=    \int_0^1[ \tau(1- \tau)]^{ \lam}
g(\Phi(\rho,u,s,\tau)) d \tau  \\\\
+ \theta \rho^{\theta}  e^{s}  \int_0^1[ \tau(1- \tau)]^{ \lam}
g'(\Phi(\rho,u,s,\tau))(1- 2 \tau) d \tau, \\\\

\eta_{u}= \rho \int_0^1[ \tau(1- \tau)]^{ \lam}
g'(\Phi(\rho,u,s,\tau)) d \tau,  \\\\

\eta_{s}= \rho^{\theta+1}  e^{s}  \int_0^1[ \tau(1- \tau)]^{ \lam}
g'(\Phi(\rho,u,s,\tau))(1- 2 \tau) d \tau,
\end{array}\right.
\end{equation}
\begin{equation}\left\{
\label{2.58}
\begin{array}{l}
\eta_{\rho \rho}= (\theta+\theta^{2}) \rho^{\theta-1}  e^{s}  \int_0^1[ \tau(1- \tau)]^{ \lam}
g'(\Phi(\rho,u,s,\tau))(1- 2 \tau) d \tau  \\\\
+ \theta^{2} \rho^{2 \theta-1}  e^{2s}  \int_0^1[ \tau(1- \tau)]^{ \lam}
g''(\Phi(\rho,u,s,\tau))(1- 2 \tau)^{2} d \tau, \\\\

\eta_{\rho u}=  \int_0^1[ \tau(1- \tau)]^{ \lam}
g'(\Phi(\rho,u,s,\tau)) d \tau  \\\\
+ \theta \rho^{\theta}  e^{s}  \int_0^1[ \tau(1- \tau)]^{ \lam}
g''(\Phi(\rho,u,s,\tau))(1- 2 \tau) d \tau, \\\\

\eta_{u u}= \rho \int_0^1[ \tau(1- \tau)]^{ \lam}
g''(\Phi(\rho,u,s,\tau)) d \tau,  \\\\
\end{array}\right.
\end{equation}
and
\begin{equation}\left\{
\label{2.59}
\begin{array}{l}
\eta_{\rho s}=(1+ \theta) \rho^{\theta}  e^{s} \int_0^1[ \tau(1- \tau)]^{ \lam}
g'(\Phi(\rho,u,s,\tau)) (1- 2 \tau) d \tau  \\\\
+ \theta \rho^{2 \theta}  e^{2s}  \int_0^1[ \tau(1- \tau)]^{ \lam}
g''(\Phi(\rho,u,s,\tau))(1- 2 \tau)^{2} d \tau, \\\\

\eta_{ss}= \rho^{1+ \theta}  e^{s} \int_0^1[ \tau(1- \tau)]^{ \lam}
g'(\Phi(\rho,u,s,\tau)) (1- 2 \tau) d \tau  \\\\
+  \rho^{2 \theta+1}  e^{2s}  \int_0^1[ \tau(1- \tau)]^{ \lam}
g''(\Phi(\rho,u,s,\tau))(1- 2 \tau)^{2} d \tau, \\\\

\eta_{us}= \rho^{1+ \theta}  e^{s} \int_0^1[ \tau(1- \tau)]^{ \lam}
g''(\Phi(\rho,u,s,\tau)) (1- 2 \tau) d \tau.
\end{array}\right.
\end{equation}
Multiplying the first equation in (\ref{1.17}) by $\eta_{\rho}$, the equation (\ref{1.42}) by $\eta_{m}$, and (\ref{1.21}) by $\eta_{s}$, then adding the result, we have (for simplicity, we omit the superscripts $\E $ and $ \delta$)
\begin{equation}\begin{array}{ll}
\label{2.60}
 \eta( \rho,u,s)_{t}+
 (q_{1}(\rho,u,s)
 - 2 \delta q_{2}(\rho,u,s))_{x} \\\\
 - (q_{1}(\rho,u,s)
 - 2 \delta q_{2}(\rho,u,s))_{s} s_{x}
 + \F{(\rho-2 \delta)}{\rho}u s_{x} \eta_{s}( \rho,u,s)
 \\\\
 = \E \eta( \rho,u,s)_{xx}+ 2 \E \F{\rho_{x}}{\rho} s_{x} \eta_{s}
 + 2 \E \F{\rho_{x}}{\rho} u_{x} \eta_{u}
 \\\\- \E ( \rho,u,s) \cdot \nabla^{2} \eta(\rho,u,s) \cdot (\rho,u,s)^{T},
\end{array}
\end{equation}
where
\begin{equation}
 \label{2.61}
\nabla^{2} \eta(\rho,u,s) =\left( \begin{array}{ccc}
\eta_{\rho \rho} & \eta_{\rho u} & \eta_{\rho s} \\\\
 \eta_{\rho u} & \eta_{u u} & \eta_{u s} \\\\
\eta_{\rho s} & \eta_{u s} & \eta_{ss}
\end{array} \right).
\end{equation}
Let
\begin{equation}
 \label{2.62}
 2 \E \F{\rho_{x}}{\rho} s_{x} \eta_{s}
 + 2 \E \F{\rho_{x}}{\rho} u_{x} \eta_{u}
 - \E ( \rho,u,s) \cdot \nabla^{2} \eta(\rho,u,s) \cdot (\rho,u,s)^{T}
 =I_{1}+I_{2}
 \end{equation}
 where $I_{1}, I_{2}$ be the sets of all functions appeared in
 (\ref{2.57})-(\ref{2.59}) with $g'(\Phi) $ and $g''(\Phi)$ respectively. Then
\begin{equation}\begin{array}{ll}
 \label{2.63}
 I_{1}= 2 \E \F{\rho_{x}}{\rho} s_{x} \eta_{s}+ 2 \E \F{\rho_{x}}{\rho} u_{x} \eta_{u} \\\\
 - \E \Big((\theta+\theta^{2}) \rho^{\theta-1}  e^{s} \rho^{2}_{x} +
 2 (1+ \theta) \rho^{\theta}  e^{s} \rho_{x}s_{x}
 + \rho^{1+ \theta}  e^{s} s^{2}_{x}  \Big) \\\\
\cdot \int_0^1[ \tau(1- \tau)]^{ \lam}
g'(\Phi(\rho,u,s,\tau))(1- 2 \tau) d \tau \\\\
- 2 \E \rho_{x} u_{x}  \int_0^1[ \tau(1- \tau)]^{ \lam}
g'(\Phi(\rho,u,s,\tau)) d \tau \\\\
= - \E \Big((\theta+\theta^{2}) \rho^{\theta-1}  e^{s} \rho^{2}_{x} +
 2  \theta \rho^{\theta}  e^{s} \rho_{x}s_{x}
 + \rho^{1+ \theta}  e^{s} s^{2}_{x}  \Big) \\\\
 \cdot \int_0^1[ \tau(1- \tau)]^{ \lam}
g'(\Phi(\rho,u,s,\tau))(1- 2 \tau) d \tau \\\\
=  - \E \Big((\theta+\theta^{2}) \rho^{\theta-1}  e^{s} \rho^{2}_{x} +
 2  \theta \rho^{\theta}  e^{s} \rho_{x}s_{x}
 + \rho^{1+ \theta}  e^{s} s^{2}_{x}  \Big) \\\\
 \cdot \F{2}{\lam+1} \rho^{\theta}  e^{s} \int_0^1[ \tau(1- \tau)]^{ \lam+1}
g''(\Phi(\rho,u,s,\tau)) d \tau \\\\

=  - \E \F{2}{\lam+1} \Big((\theta+\theta^{2}) \rho^{2 \theta-1}  e^{2s} \rho^{2}_{x} +
 2  \theta \rho^{2 \theta}  e^{2s} \rho_{x}s_{x}
 + \rho^{1+2 \theta}  e^{2s} s^{2}_{x}  \Big) \\\\
 \cdot   \int_0^1[ \tau(1- \tau)]^{ \lam+1}
g''(\Phi(\rho,u,s,\tau)) d \tau
\end{array}
 \end{equation}
 because
 \begin{equation}\begin{array}{ll}
\label{2.64}
    \int_0^1[ \tau(1- \tau)]^{ \lam}
g'(\Phi(\rho,u,s,\tau))(1- 2 \tau) d \tau \\\\
=  \F{1}{\lam+1}  \int_0^1
g'(\Phi(\rho,u,s,\tau)) d ( [ \tau(1- \tau)]^{ \lam+1}) \\\\
= \F{2}{\lam+1} \rho^{\theta}  e^{s} \int_0^1[ \tau(1- \tau)]^{ \lam+1}
g''(\Phi(\rho,u,s,\tau)) d \tau,
\end{array}
\end{equation}
 and
\begin{equation}\begin{array}{ll}
\label{2.65}
I_{2}= - \E \Big( \theta^{2} \rho^{2 \theta-1}  e^{2s} \rho^{2}_{x}
+2 \theta \rho^{2 \theta}  e^{2s} \rho_{x} s_{x}
+  \rho^{2 \theta+1}  e^{2s} s^{2}_{x} \Big) \\\\
  \cdot \int_0^1[ \tau(1- \tau)]^{ \lam}
g''(\Phi(\rho,u,s,\tau))(1- 2 \tau)^{2} d \tau  \\\\

- 2 \E \theta \rho^{\theta}  e^{s}  \int_0^1[ \tau(1- \tau)]^{ \lam}
g''(\Phi(\rho,u,s,\tau))(1- 2 \tau) d \tau  \rho_{x} u_{x} \\\\

- \E \rho \int_0^1[ \tau(1- \tau)]^{ \lam}
g''(\Phi(\rho,u,s,\tau)) d \tau   u^{2}_{x} \\\\

- 2 \E \rho^{1+ \theta}  e^{s} \int_0^1[ \tau(1- \tau)]^{ \lam}
g''(\Phi(\rho,u,s,\tau)) (1- 2 \tau) d \tau   u_{x} s_{x}.
\end{array}
\end{equation}
Since
\begin{equation}
\label{2.66}
      \E \rho^{\theta} |\rho_{x} u_{x}| \leq
      \E \rho^{\gamma-2} \rho^{2}_{x}+ \E \rho u^{2}_{x}, \quad
      \E \rho^{1+ \theta} |u_{x} s_{x}| \leq
      \E \rho^{\gamma} s^{2}_{x}+ \E \rho u^{2}_{x},
\end{equation}
then the last three terms in (\ref{2.65})
\begin{equation}\begin{array}{ll}
\label{2.67}
E(x,t) = 2 \E \theta \rho^{\theta}  e^{s}  \int_0^1[ \tau(1- \tau)]^{ \lam}
g''(\Phi(\rho,u,s,\tau))(1- 2 \tau) d \tau  \rho_{x} u_{x} \\\\

- \E \rho \int_0^1[ \tau(1- \tau)]^{ \lam}
g''(\Phi(\rho,u,s,\tau)) d \tau   u^{2}_{x} \\\\

- 2 \E \rho^{1+ \theta}  e^{s} \int_0^1[ \tau(1- \tau)]^{ \lam}
g''(\Phi(\rho,u,s,\tau)) (1- 2 \tau) d \tau   u_{x} s_{x}
\end{array}
\end{equation}
are uniformly bounded in $L^{1}_{loc}(R \times R^{+})$ due to the estimates
in (\ref{2.25}).

Moreover, the terms on the left-hand side of (\ref{2.60})
\begin{equation}
\label{2.68}
F(x,t)= - (q_{1}(\rho,u,s)
 - 2 \delta q_{2}(\rho,u,s))_{s} s_{x}
 + \F{(\rho-2 \delta)}{\rho}u s_{x} \eta_{s}( \rho,u,s)
\end{equation}
are uniformly bounded in $L^{1}_{loc}(R \times R^{+})$.

Thus we have from (\ref{2.60}),(\ref{2.63})-(\ref{2.64}) and (\ref{2.67})-(\ref{2.68}) that
\begin{equation}\begin{array}{ll}
 \label{2.69}
    \E \F{2}{\lam+1} \Big((\theta+\theta^{2}) \rho^{2 \theta-1}  e^{2s} \rho^{2}_{x} +
 2  \theta \rho^{2 \theta}  e^{2s} \rho_{x}s_{x}
 + \rho^{1+2 \theta}  e^{2s} s^{2}_{x}  \Big) \\\\
 \cdot   \int_0^1[ \tau(1- \tau)]^{ \lam+1}
g''(\Phi(\rho,u,s,\tau)) d \tau  \\\\
+ \E \Big( \theta^{2} \rho^{2 \theta-1}  e^{2s} \rho^{2}_{x}
+2 \theta \rho^{2 \theta}  e^{2s} \rho_{x} s_{x}
+  \rho^{2 \theta+1}  e^{2s} s^{2}_{x} \Big) \\\\
  \cdot \int_0^1[ \tau(1- \tau)]^{ \lam}
g''(\Phi(\rho,u,s,\tau))(1- 2 \tau)^{2} d \tau \\\\
= \E \eta(\rho,u,s)_{xx} + E(x,t)-F(x,t) \\\\-
\eta( \rho,u,s)_{t}- (q_{1}(\rho,u,s)- 2 \delta q_{2}(\rho,u,s))_{x}.
\end{array}
 \end{equation}
 Multiplying the right-hand side (we write it as $R(x,t)$ ) of Equation (\ref{2.69}) by $ \phi$, where $\phi$ is given in (\ref{2.24}), and integrating over $R \times R^+ $, we may obtain
\begin{equation}\begin{array}{ll}
\label{2.70}
|\int_{0}^{\infty} \int_{-\infty}^{\infty} R(x,t)  \phi(x,t) dxdt|
=| \int_{0}^{\infty} \int_{-\infty}^{\infty} \eta(\rho,u,s) \phi(x,t)_{xx} \\\\
+ (q_{1}(\rho,u,s)- 2 \delta q_{2}(\rho,u,s)) \phi(x,t)_{x}
+(E(x,t)-F(x,t)) \phi(x,t) dxdt|
\leq M( \phi).
\end{array}
\end{equation}
Then, if we choose $ g $, on the left-hand side of (\ref{2.69}), to be strictly convex, $g''(\Phi) \geq c >0 $
for a constant $c$, we may obtain from (\ref{2.69}) and (\ref{2.70}) that
 \begin{equation}\begin{array}{ll}
\label{2.71}
\E \Big((\theta+\theta^{2}) \rho^{2 \theta-1}   \rho^{2}_{x} +
 2  \theta \rho^{2 \theta}   \rho_{x}s_{x}
 + \rho^{1+ 2 \theta}   s^{2}_{x} \Big) \ \mbox{ are bounded in} \
L^{1}_{loc}(R \times R^{+}),
\end{array}
\end{equation}
which deduce that, for any smooth function $f$, $I_{1}+I_{2}$ in (\ref{2.62}) are
bounded in $L^{1}_{loc}(R \times R^{+})$.

Furthermore, for any $\varphi\in H^{1}_{0}(R \times R^{+})$, we have
\[  |\int_{0}^{\infty} \int_{- \infty}^{\infty} 2 \delta q_{2x} \varphi dx dt|
= |\int_{0}^{\infty} \int_{- \infty}^{\infty} 2 \delta q_{2} \varphi_{x} dx dt|
\rightarrow 0 \]
as $ \delta $ goes to zero, and from the estimates in (\ref{2.25}) that
\begin{equation}\begin{array}{ll}
\label{2.72}
|\int_{0}^{\infty} \int_{- \infty}^{\infty} \E \eta_{xx} \varphi dx dt|
= |\int_{0}^{\infty} \int_{- \infty}^{\infty}
\E (\eta_{\rho} \rho_{x}+ \eta_{u} u_{x}+ \eta_{s} s_{x}) \varphi_{x} dx dt| \\\\
\leq M (\int_{0}^{\infty} \int_{- \infty}^{\infty} \E \rho^{\gamma-2} \rho^{2}_{x} |\varphi_{x}| dx dt)^{\F{1}{2}} (\int_{0}^{\infty} \int_{- \infty}^{\infty} \E \rho^{2-\gamma} |\varphi_{x}| dx dt)^{\F{1}{2}} \\\\

+ M (\int_{0}^{\infty} \int_{- \infty}^{\infty} \E \rho u^{2}_{x} |\varphi_{x}| dx dt)^{\F{1}{2}} (\int_{0}^{\infty} \int_{- \infty}^{\infty} \E \rho^{-1} |\varphi_{x}| dx dt)^{\F{1}{2}} \\\\

+ M (\int_{0}^{\infty} \int_{- \infty}^{\infty} \E \rho^{\gamma} s^{2}_{x} |\varphi_{x}| dx dt)^{\F{1}{2}} (\int_{0}^{\infty} \int_{- \infty}^{\infty} \E \rho^{-\gamma} |\varphi_{x}| dx dt)^{\F{1}{2}} \rightarrow 0
\end{array}
\end{equation}
if we let $ \E $ go to zero much faster that $\delta$. Therefore we obtain the
proof of Lemma 2.5.

If we apply the div-curl lemma to any two pairs of weak entropy-entropy
flux given in Lemma 2.5
\begin{equation}
\begin{array}{ll}
 \label{2.73}
\Big( \eta^{(i)}(\rho^{\E,\delta},u^{\E,\delta},s^{\E,\delta}), \quad q_{1}^{(i)}(\rho^{\E,\delta},u^{\E,\delta},s^{\E,\delta}) \Big),
\end{array}
\end{equation}
where $ i=1,2 $ corresponds to $g_{i}(u)$ in (\ref{2.42}),
we have the following weak limit equations \cite{5}
\begin{equation}
\begin{array}{ll}
 \label{2.74}
\overline{ \eta^{(1)}(\rho^{\E,\delta},u^{\E,\delta},s^{\E,\delta})
\cdot q_{1}^{(2)}(\rho^{\E,\delta},u^{\E,\delta},s^{\E,\delta})

-\eta^{(2)}(\rho^{\E,\delta},u^{\E,\delta},s^{\E,\delta})
\cdot q_{1}^{(1)}(\rho^{\E,\delta},u^{\E,\delta},s^{\E,\delta})}\\\\
=\overline{ \eta^{(1)}(\rho^{\E,\delta},u^{\E,\delta},s^{\E,\delta})}
\cdot \overline{ q_{1}^{(2)}(\rho^{\E,\delta},u^{\E,\delta},s^{\E,\delta})}

-\overline{ \eta^{(2)}(\rho^{\E,\delta},u^{\E,\delta},s^{\E,\delta})}
\cdot \overline{ q_{1}^{(1)}(\rho^{\E,\delta},u^{\E,\delta},s^{\E,\delta})}.
\end{array}
\end{equation}
Using the conclusion in Lemma 2.3, we may replace $ s^{\E,\delta} $
 in (\ref{2.74}) by $s$ and have the following
weak limit equations
\begin{lemma}
 \begin{equation}
\begin{array}{ll}
 \label{2.75}
\overline{ \eta^{(1)}(\rho^{\E,\delta},u^{\E,\delta},s)
\cdot q_{1}^{(2)}(\rho^{\E,\delta},u^{\E,\delta},s)

-\eta^{(2)}(\rho^{\E,\delta},u^{\E,\delta},s)
\cdot q_{1}^{(1)}(\rho^{\E,\delta},u^{\E,\delta},s)}\\\\
=\overline{ \eta^{(1)}(\rho^{\E,\delta},u^{\E,\delta},s)}
\cdot \overline{ q_{1}^{(2)}(\rho^{\E,\delta},u^{\E,\delta},s)}

-\overline{ \eta^{(2)}(\rho^{\E,\delta},u^{\E,\delta},s)}
\cdot \overline{ q_{1}^{(1)}(\rho^{\E,\delta},u^{\E,\delta},s)},
\end{array}
\end{equation}
where $s$ is the weak-star limit of $ s^{\E,\delta} $.
\end{lemma}
{\bf Proof of Lemma 2.6.} First, using the estimate (\ref{2.33}), we have also
\begin{equation}
 \label{2.76}
\overline{ \rho^{\E,\delta} |s^{\E,\delta}-s|}=0.
\end{equation}
Second, by the entropy-entropy flux formulas given in (\ref{2.44}) and (\ref{2.45}), we have
\begin{equation}
 \label{2.77}
 |f(\rho^{\E,\delta},u^{\E,\delta},s^{\E,\delta})-
 f(\rho^{\E,\delta},u^{\E,\delta},s)| \leq M (\rho^{\E,\delta})^{\theta+1}
 |s^{\E,\delta}-s| \leq M_{1} \rho^{\E,\delta} |s^{\E,\delta}-s|,
\end{equation}
and hence the proof of Lemma 2.6, where $f$ is any one of the weak
entropies $\eta^{(i)}$ and the weak entropy fluxes  $q_{1}^{(i)}$.

{\bf Proof of (II) in Theorem 1.1.} Paying attention to the special structure of the weak entropy-entropy flux formulas given in (\ref{2.44}) and (\ref{2.45}), and letting $\omega^{\E,\delta}= \rho^{\E,\delta} e^{\F{1}{\theta} s}$, we have the following  weak limit equations from (\ref{2.75}) and Lemma 2.6,
\begin{equation}
\begin{array}{ll}
 \label{2.78}
\overline{ \eta^{(1)}(\omega^{\E,\delta},u^{\E,\delta})
\cdot q^{(2)}(\omega^{\E,\delta},u^{\E,\delta})

-\eta^{(2)}(\omega^{\E,\delta},u^{\E,\delta})
\cdot q^{(1)}(\omega^{\E,\delta},u^{\E,\delta})}\\\\
=\overline{ \eta^{(1)}(\omega^{\E,\delta},u^{\E,\delta})}
\cdot \overline{ q^{(2)}(\omega^{\E,\delta},u^{\E,\delta})}

-\overline{ \eta^{(2)}(\omega^{\E,\delta},u^{\E,\delta})}
\cdot \overline{ q^{(1)}(\omega^{\E,\delta},u^{\E,\delta})},
\end{array}
\end{equation}
where
\begin{equation}\begin{array}{ll}
\label{2.79}
\eta(\omega,u)&= \displaystyle \int_{R} g( \xi) G(\omega, \xi-u) d \xi \\\\
&= \omega \displaystyle \int_0^1[ \tau(1- \tau)]^{ \lam}
g(u+ \omega^{ \theta}- 2 \omega^{ \theta} \tau) d \tau
\end{array}
\end{equation}
and
\begin{equation}\begin{array}{ll}
\label{2.80}
q(\omega,u) = \omega \int_0^1[ \tau(1- \tau)]^{ \lam}
g(u+ \omega^{ \theta} - 2 \omega^{ \theta} \tau)(u+ \theta(1-2 \tau)
 \omega^{ \theta}) d \tau
\end{array}
\end{equation}
is a pair of weak entropy-entropy flux of the following  isentropic
gas dynamics system
\begin{equation}\left\{
\label{2.81}
\begin{array}{l}
      \omega_{t}+(\omega u)_{x}=0,        \\\\
       ( \omega u)_t+( \omega u^2+ \F{ \theta^{2}}{\gamma} \omega^{\gamma})_x =0.
\end{array}\right.
\end{equation}
Using the Young measure representation theorem from the compensated compactness
theory, we may select a subsequence (still labelled) $(\omega^{\E,\delta},u^{\E,\delta})$, and a family of positive
 measures $ \nu_{(x,t)} \in \mbox{M}(R^2)$, depending measurably on $ (x,t) \in
 K \subset R \times R^{+}$, such that
 \begin{equation}
\begin{array}{ll}
 \label{2.82}
\int_{K} \eta^{(1)}( \lam) q^{(2)}( \lam)- \eta^{(1)}( \lam) q^{(2)}( \lam)
 d \nu_{(x,t)}( \lam) \\\\
= \int_{K} \eta^{(1)}( \lam)d \nu_{(x,t)}( \lam) \cdot
\int_{K} q^{(2)}( \lam)d \nu_{(x,t)}( \lam) \\\\
- \int_{K} \eta^{(2)}( \lam)d \nu_{(x,t)}( \lam)\cdot
\int_{K} q^{(1)}( \lam)d \nu_{(x,t)}( \lam).
\end{array}
\end{equation}
With the help of the measure equations (\ref{2.82}) and
the results given in \cite{7,8,9,10,11}, we may deduce
that, for any fixed point $(x,t) \in R \times R^{+}$, the Young measure $ \nu_{(x,t)} $ is either wholly contained in the line $\omega(x,t)=0$ or
concentrated in one point $(\omega_{0}(x,t),u_{0}(x,t))$, and hence,
$ (\omega^{\E,\delta}(x,t),u^{\E,\delta}(x,t)) \rightarrow (\omega(x,t),u(x,t)) $
almost everywhere on the set $\omega_{+}= \{(x,t): \omega(x,t) >0 \}$, where
$\omega(x,t)$ is the weak-star limit of $ \omega^{\E,\delta}(x,t)$. Then,
$ (\rho^{\E,\delta}(x,t),u^{\E,\delta}(x,t),s^{\E,\delta}(x,t)) \rightarrow (\rho(x,t),u(x,t),s(x,t)) $ almost everywhere on the set $\rho_{+}= \{(x,t): \rho(x,t) >0 \}$ since $\omega^{\E,\delta}= \rho^{\E,\delta} e^{\F{1}{\theta} s}$ and $e^{\F{1}{\theta} s} > 0$.

Since the variables $(\rho,\rho u, \rho s) $ in (\ref{1.3}),
and the corresponding fluxes $(\rho u, \rho u^2+ \F{ \theta^{2}}{\gamma} \rho^{\gamma} e^{2s}, \rho u)$ are all zero at the line $\rho=0$,  we may prove that the set of functions $(\rho,u,s)$ satisfies (\ref{1.28}) and (\ref{1.29}) by letting $ \E, \delta $ in (\ref{1.17}) go to zero. Thus we complete the proof of Theorem 1.1.

\noindent {\bf Acknowledgments:}  This work was started when the author visited Heidelberg University, Germany as a Humboldt fellow and completed with the support
 of a Humboldt renewed research fellowship in University of Wurzburg, Germany. The author is very grateful to the colleagues in these two universities for their warm hospitality.



\begin{thebibliography}{01}
\addcontentsline{toc}{chapter}{Bibliography}




\bibitem[1]{1} F. Bereux, E. Bonnetier and P. Lefloch, {\em Gas dynamics system: two special cases}, SIAM J. Math. Anal. {\bf 28} (1997), 499-515.

 \bibitem[2]{2} T. Nishida, {\em Global solution for an initial-boundary-value
problem of a quasilinear hyperbolic system},
Proc. Jap. Acad., {\bf 44} (1968), 642-646.

\bibitem[3]{3} T. Nishida and J. Smoller, {\em Solutions in the large
for some nonlinear hyperbolic conservation laws},
Comm. Pure Appl. Math., {\bf 26} (1973), 183-200.

\bibitem[4]{4} J. Glimm, {\em Solutions in the large for nonlinear hyperbolic
systems of equations}, Comm. Pure Appl. Math., {\bf 18} (1965),
95-105.

\bibitem[5]{5} T. Tartar, {\em Compensated compactness and applications to
partial differential equations}, In: Research Notes in Mathematics,
Nonlinear Analysis and Mechanics, Heriot-Watt symposium, Vol. {\bf
4}, ed. R.~J. Knops, Pitman Press, London, 1979.

\bibitem[6]{6} F. Murat, {\em Compacit\'e par compensation},
 Ann. Scuola Norm. Sup. Pisa, {\bf 5} (1978), 489-507.


\bibitem[7]{7} R.~J. DiPerna, {\em Convergence of the viscosity method for isentropic gas dynamics},
Commun. Math. Phys., {\bf 91} (1983), 1-30.


\bibitem[8]{8} X.-X. Ding, G.-Q. Chen and P.-Z. Luo, {\em Convergence of the Lax-Friedrichs schemes for the isentropic gas dynamics I-II}, Acta Math. Sci., {\bf 5} (1985), 415-432, 433-472.

\bibitem[9]{9} G.-Q. Chen, {\em Convergence of the Lax-Friedrichs
scheme for isentropic gas dynamics}, Acta Math. Sci., {\bf 6} (1986), 75-120.

\bibitem[10]{10} P.~L. Lions, B. Perthame and E. Tadmor,
{\em Kinetic formulation of the isentropic gas dynamics and
p-system}, Commun. Math. Phys., {\bf 163} (1994), 415-431.

\bibitem[11]{11} P.~L. Lions, B. Perthame and P.~E. Souganidis,
{\em Existence and stability of entropy solutions for the hyperbolic
systems of isentropic gas dynamics in Eulerian and Lagrangian
coordinates}, Comm. Pure Appl. Math., {\bf 49} (1996), 599-638.

\bibitem[12]{12} Y.-G. Lu, {\em Existence of Global Entropy Solutions to a Nonstrictly  Hyperbolic System},
Arch. Rat. Mech. Anal., {\bf 178}(2005), 287-299.

\bibitem[13]{13} F.M. Huang and Z. Wang, {\em Convergence of Viscosity Solutions for Isentropic
Gas Dynamics}, SIAM J. Math. Anal., {\bf 34} (2003), 595-610.

\bibitem[14]{14} T.-P. Liu, {\em Solutions in the large for the equations of nonisentropic gas dynamics}, Indiana Univ. Math. J., {\bf 26}(1977), 797-838.

\bibitem[15]{15} B. Temple, {\em Systems of conservation laws with invariant
submanifolds}, Trans. of Am. Math. Soc., {\bf 280} (1983), 781-795.


\bibitem[16]{16} K.~N. Chueh, C.~C. Conley and J.~A. Smoller,
{\em Positive invariant regions for systems of nonlinear diffusion equations},
Indiana Univ. Math. J., {\bf 26} (1977), 372-411.

\bibitem[17]{17} S. Benzoni-Gavage and D. Serre, {\em Compacit\'e par compensation pour une classe
de syst\'emes hyperboliques de $p \geq 3 $ lois de conservation},  Rev. Mat. Iberoamericana,
{\bf 10} (1994), no. 3, 557-579.

\bibitem[18]{18} R.~J. DiPerna, {\em Convergence of approximate solutions to conservation laws},
Arch. Rat. Mech. Anal., {\bf 82} (1983), 27-70.

\bibitem[19]{19} Y.-J. Peng, {\em Solutions faibles globales pour un modele
decoulements diphasiques}, Ann. Scuola
Norm. Sup. Pisa Cl. Sci., {\bf 21} (1994), 523-540.

\bibitem[20]{20} Changjiang Zhu, {\em Global smooth solution of the
nonisentropic gas dynamics system}, Proc. Royal Soc. Edinburgh, {\bf 126A} (1996),
769-775.

\bibitem[21]{21} H. Frid, H. Holden and K.H. Karlsen, {\em $ L^{\infty} $ solutions for a model of polytropic gas
flow with diffusive entropy}, SIAM J. Math. Anal., {\bf 43} (2011), 2253-2274.

\bibitem[22]{22} Y.-G. Lu, {\em Some Results on General System of Isentropic Gas
Dynamics}, Differential Equations, {\bf 43} (2007), 130-138.

\bibitem[23]{23}  Y.-G. Lu, {\em Global Existence of Resonant Isentropic Gas Dynamics},
 Nonlinear Analysis, Real World Applications, {\bf 12}(2011), 2802-2810.

\bibitem[24]{24} Y.-G. Lu, {\em Existence of Global Bounded Weak Solutionsto a Non-Symmetric System of Keyfitz-Kranzer type},
J. Funct. Anal., { \bf 261}(2011), 2797-2815.

\bibitem[25]{25} Y.-G. Lu, {\em Existence of Global Bounded Weak Solutions
to a Non-Symmetric System of Keyfitz-Kranzer type}, J. Funct.
Anal., 264(2013), 2457-2468.

\bibitem[26]{26} D. Serre, {\em Solutions \`a variations born\'ees pour certains syst\`emes
hyperboliques de lois de conservation},  J. Diff. Eqs., {\bf 68}
(1987), 137-168.

\bibitem[27]{27} F. James, Y.-J. Peng and B. Perthame, {\em Kinetic formulation for chromatography and some other hyperbolic systems}, J. Math. Pure
Appl., {\bf 74} (1995), 367-385.

\bibitem[28]{28} Y.-G. Lu, {\em Hyperbolic Conservation Laws and the
Compensated Compactness Method}, Vol. {\bf 128}, Chapman and Hall,
CRC Press, New York, 2002.

\end{thebibliography}
\end{document}